\documentclass[11pt]{amsart}

\usepackage{multirow}
\usepackage{graphicx}
\usepackage[top=26mm, bottom=26mm, left=23mm, right=25mm]{geometry}
\usepackage[dvips]{color}
\usepackage{bm}
\usepackage{marvosym}
\theoremstyle{definition}
\newtheorem{remark}{Remark}
\numberwithin{remark}{section}

\numberwithin{definition}{section}

\newcommand{\ngauss}{q}
\newcommand{\ncheb}{p}

\newcommand{\pvct}[1]{\bm{#1}}
\newcommand{\vct}[1]{\bm{\mathsf{#1}}}

\newcommand{\pxx}{\pvct{x}}
\newcommand{\pyy}{\pvct{y}}
\newcommand{\pzz}{\pvct{z}}
\newcommand{\uu}{\vct{u}}
\newcommand{\vv}{\vct{v}}
\newcommand{\ww}{\vct{w}}

\renewcommand{\AA}{\mtx{A}}

\newcommand{\TT}{\mtx{T}}

\newcommand{\mtx}[1]{\bm{\mathsf{#1}}}

\newcommand{\mtwo}[4]{\left[\begin{array}{cc} #1 & #2 \\ #3 & #4 \end{array}\right]}
\newcommand{\vtwo}[2]{\left[\begin{array}{cc} #1 \\ #2 \end{array}\right]}

\newcommand{\pgnotate}[1]{}

\newcommand{\lsp}{\vspace{3mm}}

\begin{document}

\begin{center}
\textbf{\large A direct solver with $O(N)$ complexity for variable coefficient elliptic PDEs\\
discretized via a high-order composite spectral collocation method}

\lsp

{\small\it
\begin{tabular}{ccc}
A. Gillman &\mbox{}\hspace{5mm}\mbox{}& P.G. Martinsson\\
Department of Mathematics&&Department of Applied Mathematics\\
Dartmouth College &&University of Colorado at Boulder
\end{tabular}}

\lsp

\begin{minipage}{135mm}
\noindent\textbf{Abstract:}
A numerical method for solving elliptic PDEs with variable coefficients
on two-dimensional domains is presented. The method is based on
high-order composite spectral approximations and is designed for
problems with smooth solutions.  The resulting system of
linear equations is solved using a direct (as opposed to iterative)
solver that has optimal $O(N)$ complexity
for all stages of the computation when applied to problems with non-oscillatory
solutions such as the Laplace and the Stokes equations.
Numerical examples demonstrate that the scheme is capable of computing
solutions with relative accuracy of $10^{-10}$ or better, even for challenging
problems such as highly oscillatory Helmholtz problems and convection-dominated
convection diffusion equations.
In terms of speed, it is demonstrated that a problem with a non-oscillatory solution
that was discretized using $10^{8}$ nodes was solved in $115$ minutes on a
personal work-station with two quad-core 3.3GHz CPUs. Since the solver is direct,
and the ``solution operator'' fits in RAM, any solves beyond the first are very fast.
In the example with $10^{8}$ unknowns, solves require only $30$ seconds.
%
\end{minipage}
\end{center}

\section{Introduction}
\label{sec:intro}

\subsection{Problem formulation}

The paper describes a numerical method with optimal $O(N)$ complexity for solving boundary value
problems of the form
\begin{equation}
\label{eq:basic}
\left\{\begin{aligned}
Au(\pxx) =&\ 0\qquad &\pxx \in \Omega,\\
u(\pxx)    =&\ f(\pxx)\qquad &\pxx \in \Gamma,
\end{aligned}\right.
\end{equation}
where $\Omega$ is a rectangle in the plane with boundary $\Gamma$,
and where $A$ is a coercive elliptic partial differential operator
\begin{multline}
\label{eq:defA}
[Au](\pxx) = -c_{11}(\pxx)[\partial_{1}^{2}u](\pxx)
-2c_{12}(\pxx)[\partial_{1}\partial_{2}u](\pxx)
-c_{22}(\pxx)[\partial_{2}^{2}u](\pxx)\\
+c_{1}(\pxx)[\partial_{1}u](\pxx)
+c_{2}(\pxx)[\partial_{2}u](\pxx)
+c(\pxx)\,u(\pxx).
\end{multline}
The methodology is based on a high order composite spectral discretization
and can be modified to handle a range of different domains, including
curved ones. For problems with smooth solutions, we demonstrate that the
method can easily produce answers with ten or more correct digits.

The proposed method is based on a \textit{direct solver} which in a single sweep
constructs an approximation to the solution operator of (\ref{eq:basic}). This
gives the solver several advantages over established linear complexity
methods based on \textit{iterative solvers} (e.g.~GMRES or multigrid),
perhaps most importantly, the new method can solve problems for which iterative
methods converge slowly or not at all. The direct solver has $O(N)$
complexity for all stages of the computation. A key feature is that once
the solution operator has been built, solves can be executed
extremely rapidly, making the scheme excel when solving
a sequence of equations with the same operator but different boundary data.

\subsection{Outline of solution procedure}
\label{sec:introsummary}
The method in this paper is comprised of three steps:
\begin{enumerate}
\item The domain is first tessellated into a hierarchical tree
of rectangular patches. For each patch on the finest level, a local ``solution
operator'' is built using a dense brute force calculation. The
``solution operator'' will be defined in Section \ref{sec:DtN}; for
now we simply note that it encodes all information about the patch
that is required to evaluate interactions with other patches.
\item The larger patches are processed in an upwards pass through the tree,
where each parent can be processed once its children have been processed.
The processing of a parent node consists of forming its solution operator
by ``gluing together'' the solution  operators of its children.
\item Once the solution operators for all patches have been computed,
a solution to the PDE can be computed via a downwards pass through
the tree. This step is typically \textit{very} fast.
\end{enumerate}



\subsection{Local solution operators}
\label{sec:DtN}
The ``local solution operators'' introduced in Section \ref{sec:introsummary}
take the form of discrete approximations to the \textit{Dirichlet-to-Neumann}, or ``DtN,'' maps.
To explain what these maps do, first observe that for a given boundary function
$f$, the BVP (\ref{eq:basic}) has a unique solution $u$ (recall that we assume $A$ to be coercive).
For $\pxx \in \Gamma$, let $g(\pxx) = u_{n}(\pxx)$ denote the normal derivative in the outwards direction
of $u$ at $\pxx$. The process for constructing the function $g$ from $f$ is linear,
we write it as
$$
g = T\,f.
$$
Or, equivalently,
$$
T\,\colon\,u|_{\Gamma} \mapsto u_{n}|_{\Gamma},
\quad \mbox{where} \ u\ \mbox{satisfies} \
Au = 0\ \mbox{in}\ \Omega.
$$

From a mathematical perspective, the map $T$ is a slightly unpleasant object;
it is a hyper-singular integral operator whose kernel exhibits complicated
behavior near the corners of $\Gamma$. A key observation is that in the
present context, these difficulties can be
ignored since we limit attention to functions that are smooth. In a
sense, we only need to accurately represent the projection of the
``true'' operator $T$ onto a space of smooth functions (that in particular
do not have any corner singularities).

Concretely, given a square box $\Omega_{\tau}$ we represent a boundary
potential $u|_{\Gamma}$ and a boundary flux $u_{n}|_{\Gamma}$ via
tabulation at a set of $r$ tabulation
points on each side. (For a leaf box, we use $r$ Gaussian nodes.)
The DtN operator $T^{\tau}$ is then represented simply as a dense
matrix $\mtx{T}^{\tau}$ of size $4r \times 4r$ that maps tabulated boundary
potentials to the corresponding tabulated boundary fluxes.

%

\subsection{Computational complexity}
A straight-forward implementation of the direct solver outlined in
Sections \ref{sec:introsummary} and \ref{sec:DtN} in which all solution operators
$\mtx{T}^{\tau}$ are treated as general dense matrices has asymptotic
complexity $O(N^{1.5})$ for the ``build stage'' where the solution
operators are constructed, and $O(N \log N)$ complexity for the
``solve stage'' where the solution operator is applied for a given
set of boundary data \cite{2012_martinsson_spectralcomposite_jcp}.
This paper demonstrates that by exploiting internal structure in
these operators, they can be stored and manipulated efficiently,
resulting in optimal $O(N)$ overall complexity.

To be precise, the internal structure exploited is that the
off-diagonal blocks of the dense solution operators can to high
accuracy be approximated by matrices of low rank. This property
is a result of the fact that for a patch $\Omega_{\tau}$, the
matrix $\mtx{T}^{\tau}$ is a discrete approximation of the continuum
DtN operator $T^{\tau}$, which is an integral operator whose kernel
is smooth away from the diagonal.

\begin{remark}
The proposed solver can with slight modifications be applied to
non-coercive problems such as the Helmholtz equation. If the equation
is kept fixed while $N$ is increased, $O(N)$ complexity is retained.
However, in the context of elliptic problems with oscillatory solutions,
it is common to scale $N$ to the wave-length so
that the number of discretization points per wave-length is fixed as
$N$ increases. Our accelerated technique will
in this situation lead to a practical speed-up, but will have the
same $O(N^{1.5})$ asymptotic scaling as the basic method that does
not use fast operator algebra.
\end{remark}


\subsection{Prior work}
The direct solver outlined in Section \ref{sec:introsummary} is
an evolution of a sequence of direct solvers for integral equations
dating back to \cite{2005_martinsson_fastdirect} and later
\cite{2009_martinsson_ACTA,m2011_1D_survey,2012_ho_greengard_fastdirect,2002_chen_direct_lippman_schwinger,2013_yu_chen_totalwave}.
The common idea is to build a global solution operator by splitting
the domain into a hierarchical tree of patches, build a local solution
operator for each ``leaf'' patch, and then build solution operators for
larger patches via a hierarchical merge procedure in a sweep over the tree
from smaller to larger patches. In the context of integral equations, the
``solution operator'' is a type of scattering matrix while in the present
context, the solution operator is a DtN operator.

The direct solvers \cite{2005_martinsson_fastdirect,2009_martinsson_ACTA,m2011_1D_survey,2012_ho_greengard_fastdirect,2002_chen_direct_lippman_schwinger},
designed for dense systems, are conceptually related
to earlier work on direct solvers for sparse systems
arising from finite difference and finite element discretizations of elliptic
PDEs such as the classical nested dissection method of George \cite{george_1973,hoffman_1973}
and the multifrontal methods by Duff and others \cite{1989_directbook_duff}.
These techniques typically require $O(N^{1.5})$ operations to construct the
LU-factorization of a sparse coefficient matrix arising from the discretization
of an elliptic PDE on a planar domain, with the
dominant cost being the formation of Schur complements and LU-factorizations
of dense matrices of size up to $O(N^{0.5}) \times O(N^{0.5})$. It was in the
last several years demonstrated that these dense matrices have internal structure
that allows the direct solver to be accelerated to linear or close to linear
complexity, see, e.g., \cite{2009_xia_superfast,Adiss,2007_leborne_HLU,2009_martinsson_FEM,2011_ying_nested_dissection_2D}.
These accelerated nested dissection methods are closely related to the fast direct
solver presented in this manuscript. An important difference is that the method in the
present paper allows high order discretizations to be used without increasing the
cost of the direct solver. To be technical, the solvers in
\cite{2009_xia_superfast,Adiss,2007_leborne_HLU,2009_martinsson_FEM,2011_ying_nested_dissection_2D}
are based on an underlying finite difference or finite element discretization.
High order discretization in this context leads to large frontal matrices
(since the ``dividers'' that partition the grid have to be wide), and consequently
very high cost of the LU-factorization.

Our discretization scheme is related to earlier work on spectral collocation methods
on composite (``multi-domain'') grids, such as, e.g., \cite{1998_kopriva_multidomain,2000_hesthaven_pseudospectral},
and in particular Pfeiffer \textit{et al} \cite{2003_pfeiffer_spectralmultidomain}.
For a detailed review of the similarities and differences, see \cite{2012_martinsson_spectralcomposite_jcp}.

An $O(N^{1.5})$ complexity version of the direct solver described in this paper
was presented in \cite{2012_martinsson_spectralcomposite_jcp} which in
turn is based on \cite{2012_martinsson_spectralcomposite}. In addition to the improvement
in complexity, this paper describes a new representation
of the local solution operators that leads to cleaner implementation of the
direct solvers and allows greater flexibility in executing the leaf computation,
see Remark \ref{re:chebvsgauss}.

\subsection{Outline of paper}
Section \ref{sec:discretization} introduces the mesh of Gaussian nodes
that forms our basic computational grid.
Sections \ref{sec:leaf}, \ref{sec:merge}, and \ref{sec:hierarchy} describe a
relatively simple direct solver with $O(N^{1.5})$ complexity.
Sections \ref{sec:HBS}, \ref{sec:fastHBS}, and \ref{sec:accel} describe how
to improve the asymptotic complexity of the direct solver from $O(N^{1.5})$
to $O(N)$ by exploiting internal structure in certain dense matrices.
Section \ref{sec:numerics} describes numerical examples and Section \ref{sec:conc}
summarizes the key findings.


\section{Discretization}
\label{sec:discretization}

Partition the domain $\Omega$ into a collection of square (or
possibly rectangular) boxes, called \textit{leaf boxes}.
On the edges of each leaf, place $\ngauss$ Gaussian interpolation
points. The size of the leaf boxes, and the parameter $\ngauss$
should be chosen so that any potential solution $u$ of
(\ref{eq:basic}), as well as its first and second derivatives,
can be accurately interpolated from their values at these points
($\ngauss=21$ is often a good choice).
Let $\{\pxx_{k}\}_{k=1}^{N}$ denote the collection of interpolation
points on all boundaries.

\begin{figure}
\includegraphics[width=\textwidth]{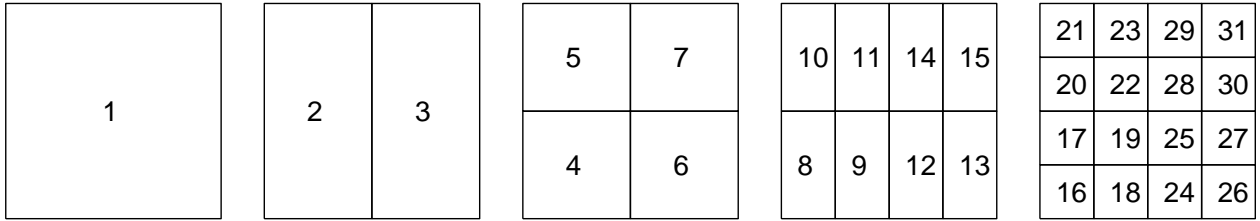}
\caption{The square domain $\Omega$ is split into $4 \times 4$ leaf boxes.
These are then gathered into a binary tree of successively larger boxes
as described in Section \ref{sec:thealgorithm}. One possible enumeration
of the boxes in the tree is shown, but note that the only restriction is
that if box $\tau$ is the parent of box $\sigma$, then $\tau < \sigma$.}
\label{fig:tree_numbering}
\end{figure}

Next construct a binary tree on the collection of leaf boxes by
hierarchically merging them, making sure that all boxes on
the same level are roughly of the same size, cf.~Figure
\ref{fig:tree_numbering}.  The boxes should be ordered so
that if $\tau$ is a parent of a box $\sigma$, then $\tau < \sigma$. We
also assume that the root of the tree (i.e.~the full box $\Omega$) has
index $\tau=1$. We let $\Omega_{\tau}$ denote the domain associated with box $\tau$.

With each box $\tau$, we define two index vectors $I_{\rm i}^{\tau}$ and $I_{\rm e}^{\tau}$ as follows:

\vspace{1mm}

\begin{itemize}
\item[$I_{\rm e}^{\tau}$] A list of all \textit{exterior} nodes of $\tau$.
In other words,
$k \in I_{\rm e}^{\tau}$ iff $\pxx_{k}$ lies on the boundary of $\Omega_{\tau}$.

\vspace{1mm}

\item[$I_{\rm i}^{\tau}$]
For a parent $\tau$, $I_{\rm i}^{\tau}$ is a list of all its \textit{interior} nodes
that are not interior nodes of its children.
For a leaf $\tau$, $I_{\rm i}^{\tau}$ is empty.
\end{itemize}

\vspace{1mm}

\noindent
Let $\uu \in \mathbb{R}^{N}$ denote a vector holding approximations
to the values of $u$ of (\ref{eq:basic}), in other words,
$$
\uu(k) \approx u(\pxx_{k}).
$$
Finally, let $\vv \in \mathbb{R}^{N}$ denote a vector holding approximations to
the boundary fluxes of the solution $u$ of (\ref{eq:basic}), in other words
$$
\vv(k) \approx \left\{
\begin{aligned}
\partial_{2}u(\pxx_{k}),\qquad&\mbox{when }\pxx_{j}\mbox{ lies on a horizontal edge,}\\
\partial_{1}u(\pxx_{k}),\qquad&\mbox{when }\pxx_{j}\mbox{ lies on a vertical edge.}
\end{aligned}\right.
$$
Note the $\vv(k)$
represents an \textit{outgoing flux} on certain boxes and an \textit{incoming flux}
on others. This is a deliberate choice to avoid problems with signs when matching fluxes
of touching boxes.


\section{Constructing the Dirichlet-to-Neumann map for a leaf}
\label{sec:leaf}

This section describes a spectral method for computing a discrete approximation to the DtN
map $T^{\tau}$ associated with a leaf box $\Omega_{\tau}$. In other words, if $u$ is a
solution of (\ref{eq:basic}), we seek a matrix $\mtx{T}^{\tau}$ of size $4\ngauss \times 4\ngauss$
such that
\begin{equation}
\label{eq:Ttau}
\vv(I_{\rm e}^{\tau}) \approx \mtx{T}^{\tau}\,\uu(I_{\rm e}^{\tau}).
\end{equation}
Conceptually, we proceed as follows: Given a vector $\uu(I_{\rm e}^{\tau})$ of potential
values tabulated on the boundary of $\Omega_{\tau}$, form for each side the unique polynomial
of degree at  most $\ngauss-1$ that interpolates the $\ngauss$ specified values of $u$. This yields Dirichlet
boundary data on $\Omega_{\tau}$ in the form of four polynomials. Solve the restriction
of (\ref{eq:basic}) to $\Omega_{\tau}$ for the specified boundary data using a spectral method
on a local tensor product grid of $\ngauss \times \ngauss$ \textit{Chebyshev nodes.} The vector
$\vv(I_{\rm e}^{\tau})$ is obtained by spectral differentiation of the local solution, and
then re-tabulating the boundary fluxes to the Gaussian nodes in $\{\pxx_{k}\}_{k \in I_{\rm e}^{\tau}}$.

We give details of the construction in Section \ref{sec:const}, but as a preliminary step,
we first review a classical spectral collocation method for the local solve in Section \ref{sec:prelim}

\begin{remark}
\label{re:chebvsgauss}
Chebyshev nodes are ideal for the leaf computations, and it is in principle also possible
to use Chebyshev nodes to represent all boundary-to-boundary ``solution operators'' such
as, e.g., $\mtx{T}^{\tau}$ (indeed, this was the approach taken in the
first implementation of the proposed method \cite{2012_martinsson_spectralcomposite_jcp}).
However, there are at least two substantial benefits to using Gaussian nodes that justify
the trouble to retabulate the operators. First, the procedure for merging boundary operators
defined for neighboring boxes is much cleaner and involves less bookkeeping since the
Gaussian nodes do not include the corner nodes. (Contrast Section 4 of \cite{2012_martinsson_spectralcomposite_jcp}
with Section \ref{sec:merge}.)  Second, and more importantly, the use of the Gaussian
nodes allows for interpolation between different discretizations. Thus the method can
easily be extended to have local refinement when necessary, see Remark \ref{remark:refine}.
%
\end{remark}

\subsection{Spectral discretization}
\label{sec:prelim}
Let $\Omega_{\tau}$ denote a rectangular subset of $\Omega$ with boundary $\Gamma_{\tau}$, and consider
the local Dirichlet problem
\begin{align}
\label{eq:differential}
[Au](\pxx) =&\ 0,\qquad &\pxx \in \Omega_{\tau}\\
   u(\pxx) =&\ h(\pxx),\qquad &\pxx \in \Gamma_{\tau},
\end{align}
where the elliptic operator $A$ is defined by (\ref{eq:defA}).
We will construct an approximate solution to (\ref{eq:differential})
using a classical spectral collocation method described in, e.g.,
Trefethen \cite{2000_trefethen_spectral_matlab}:
First, pick a small integer $\ncheb$ and let
$\{\pzz_{k}\}_{k=1}^{\ncheb^{2}}$
denote the nodes in a tensor product grid of $\ncheb\times \ncheb$ Chebyshev
nodes on $\Omega_{\tau}$.
Let $\mtx{D}^{(1)}$ and $\mtx{D}^{(2)}$ denote spectral differentiation matrices corresponding to
the operators $\partial/\partial x_1$ and $\partial/\partial x_2$, respectively.
The operator (\ref{eq:defA}) is then locally approximated via the $\ncheb^{2} \times \ncheb^{2}$ matrix
\begin{equation}
\AA =
-\mtx{C}_{11}\bigl(\mtx{D}^{(1)}\bigr)^{2}
-2\mtx{C}_{12}\mtx{D}^{(1)}\mtx{D}^{(2)}
-\mtx{C}_{22}\bigl(\mtx{D}^{(2)}\bigr)^{2}
+\mtx{C}_{1}\mtx{D}^{(1)}
+\mtx{C}_{2}\mtx{D}^{(2)}
+\mtx{C},
\end{equation}
where $\mtx{C}_{11}$ is the diagonal matrix with diagonal entries $\{c_{11}(\pzz_{k})\}_{k=1}^{\ncheb^{2}}$,
and the other matrices $\mtx{C}_{ij}$, $\mtx{C}_{i}$, $\mtx{C}$ are defined analogously.

\begin{figure}
\begin{tabular}{ccc}
\includegraphics[height=50mm]{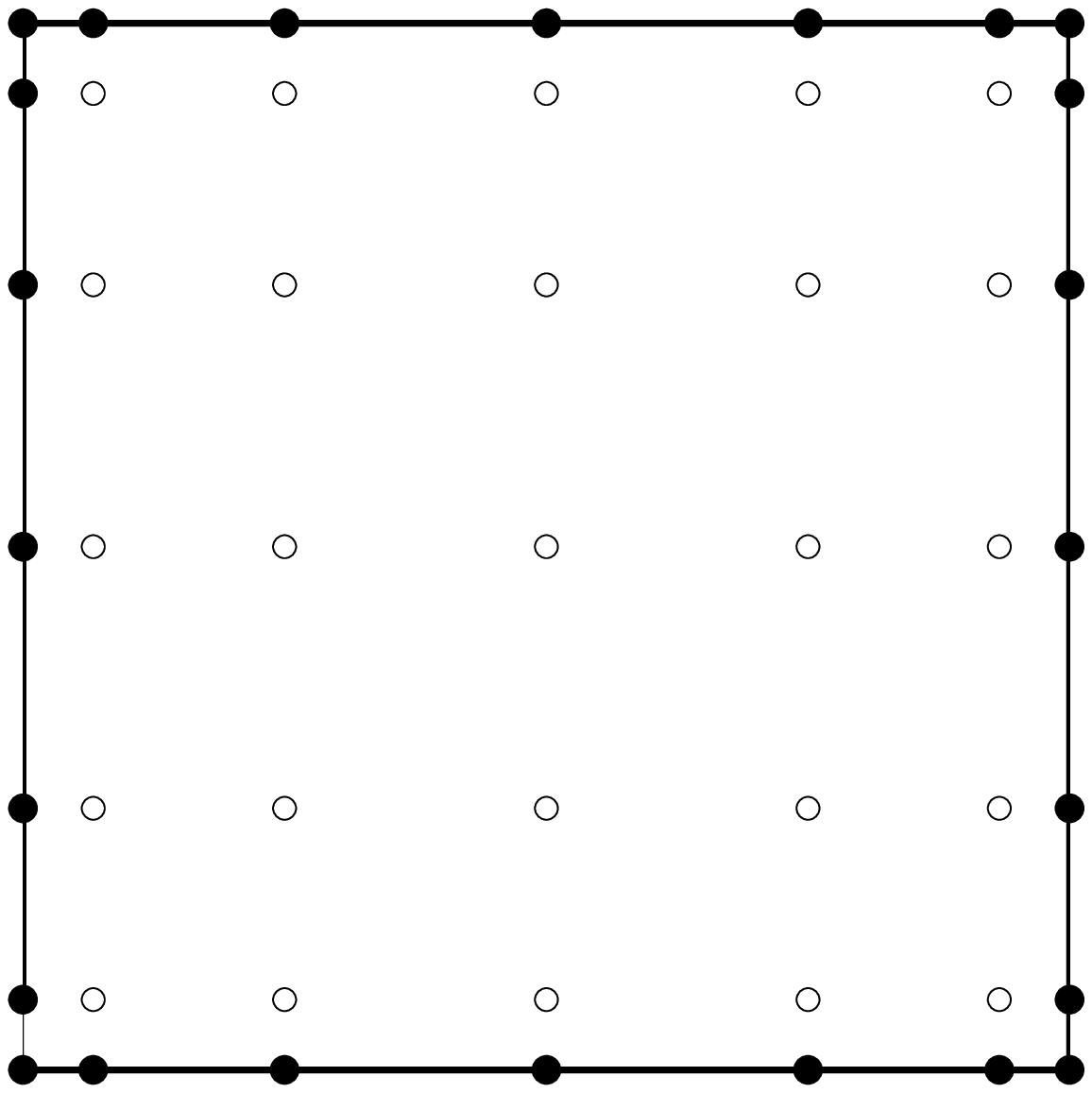}
& \mbox{}\hspace{10mm}\mbox{} &
\includegraphics[height=50mm]{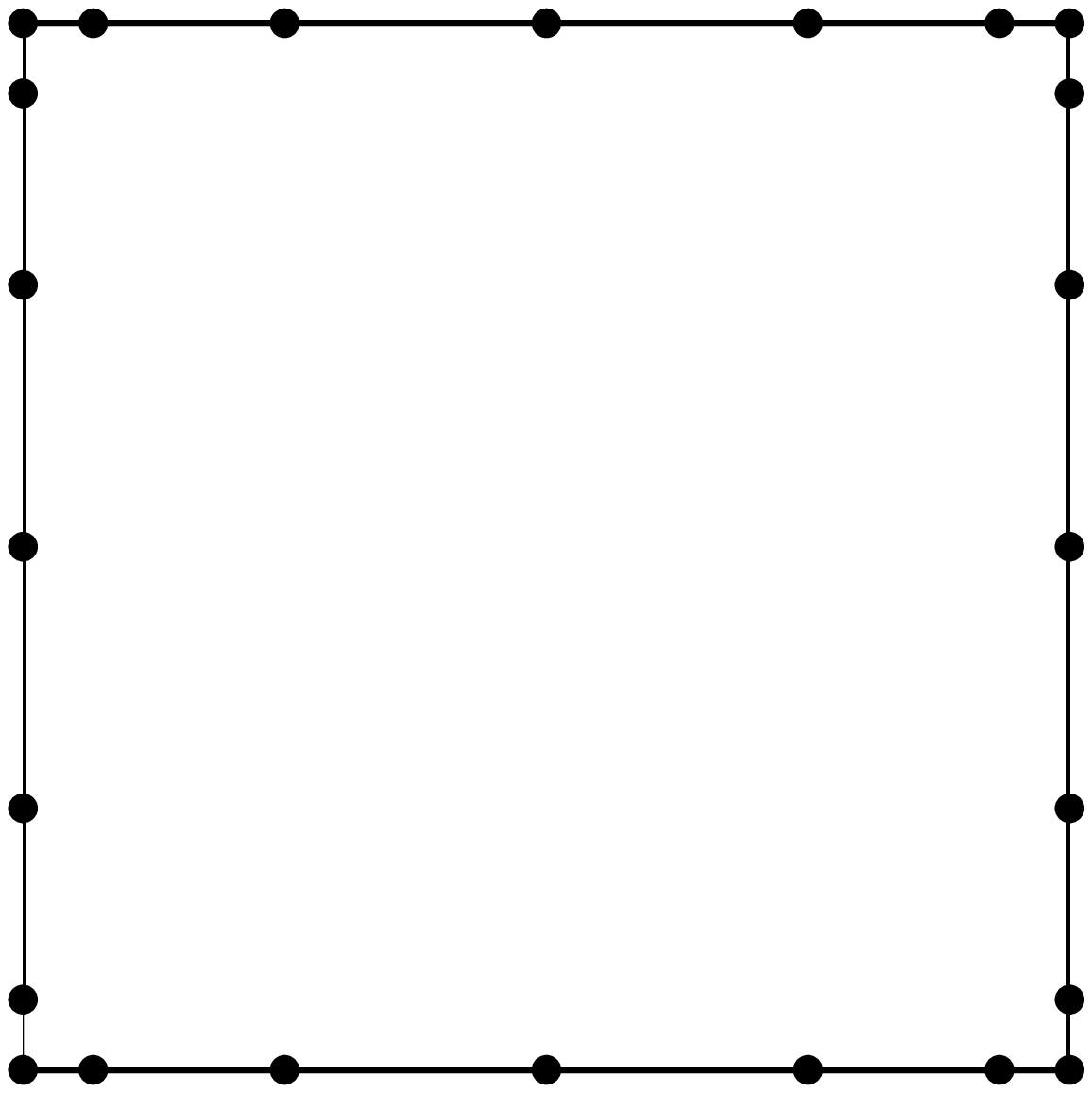}\\
(a) && (b)
\end{tabular}
\caption{Notation for the leaf computation in Section \ref{sec:leaf}.
(a) A leaf before elimination of interior (white) nodes.
(b) A leaf after elimination of interior nodes.}
\label{fig:leaf}
\end{figure}

Let $\ww \in \mathbb{R}^{\ncheb^{2}}$ denote a vector holding the desired
approximate solution of (\ref{eq:differential}). We populate all entries corresponding
to boundary nodes with the Dirichlet data from $h$, and then enforce a spectral
collocation condition at the interior nodes. To formalize, let us partition the index set
$$
\{1,\,2,\,\dots,\,\ncheb^{2}\} = J_{\rm e}\cup J_{\rm i}
$$
in such a way that $J_{\rm e}$ contains the $4(\ncheb-1)$ nodes on the boundary of $\Omega_{\tau}$, and
$J_{\rm i}$ denotes the set of $(\ncheb-2)^{2}$ interior nodes, see Figure \ref{fig:leaf}(a).
Then partition the vector $\ww$ into two parts corresponding to internal and
exterior nodes via
$$
\ww_{\rm i} = \ww(J_{\rm i}),\quad
\ww_{\rm e} = \ww(_{\rm e}).
$$J
Analogously, partition $\AA$ into four parts via
$$
\mtx{A}_{\rm i,i} = \mtx{A}(J_{\rm i},J_{\rm i}),\quad
\mtx{A}_{\rm i,e} = \mtx{A}(J_{\rm i},J_{\rm e}),\quad
\mtx{A}_{\rm e,i} = \mtx{A}(J_{\rm e},J_{\rm i}),\quad
\mtx{A}_{\rm e,e} = \mtx{A}(J_{\rm e},J_{\rm e}).
$$
The potential at the exterior nodes is now given directly from the
boundary condition:
$$
\ww_{\rm e} = \left[h(\pzz_{k})\right]_{k \in J_{\rm e}}.
$$
For the internal nodes, we enforce the PDE (\ref{eq:differential}) via direct collocation:
\begin{equation}
\label{eq:A_interior}
\AA_{\rm i,i}\,\ww_{\rm i} +
\AA_{\rm i,e}\,\ww_{\rm e} = \vct{0}.
\end{equation}
Solving (\ref{eq:A_interior}) for $\ww_{\rm i}$, we find
\begin{equation}
\label{eq:leaf_solve}
\ww_{\rm i} = -\AA_{\rm i,i}^{-1}\,\AA_{\rm i,e}\,\ww_{\rm e},
\end{equation}

\subsection{Constructing the approximate DtN}
\label{sec:const}
Now that we know how to approximately solve the local Dirichlet problem
(\ref{eq:differential}) via a local spectral method, we can
build a matrix $\mtx{T}^{\tau}$ such that (\ref{eq:Ttau})
holds to high accuracy. The starting point is a vector $\uu(I_{\tau}) \in \mathbb{R}^{4\ngauss}$ of
tabulated potential values on the boundary of $\Omega_{\tau}$. We will
construct the vector $\vv(I_{\tau}) \in \mathbb{R}^{4\ngauss}$ via four linear maps. The combination
of these maps is the matrix $\mtx{T}^{\tau}$. We henceforth assume that
the spectral order of the local Chebyshev grid matches the order of the
tabulation on the leaf boundaries so that $\ncheb=\ngauss$.

\lsp

\noindent
\textbf{\textit{Step 1 --- re-tabulation from Gaussian nodes to Chebyshev nodes:}}
For each side of $\Omega_{\tau}$, form the unique interpolating polynomial
of degree at most $\ngauss-1$ that interpolates the $\ngauss$ potential values on
that side specified by $\uu(I_{\rm e}^{\tau})$. Now evaluate these polynomials at the
boundary nodes of a $\ngauss \times \ngauss$ Chebyshev grid on $\Omega_{\tau}$.
Observe that for a corner node, we may in the general case get conflicts. For
instance, the potential at the south-west corner may get one value from extrapolation
of potential values on the south border, and one value from extrapolation of the
potential values on the west border. We resolve such conflicts by assigning the corner
node the average of the two possibly different values. (In practice, essentially no
error occurs since we know that the vector $\uu(I_{\rm e}^{\tau})$ tabulates an
underlying function that is continuous at the corner.)

\lsp

\noindent
\textbf{\textit{Step 2 --- spectral solve:}}
Step 1 populates the boundary nodes of the $\ngauss \times \ngauss$ Chebyshev grid with
Dirichlet data. Now determine the potential at all interior points on the Chebyshev
grid by executing a local spectral solve, cf.~equation (\ref{eq:leaf_solve}).

\lsp

\noindent
\textbf{\textit{Step 3 --- spectral differentiation:}}
After Step 2, the potential is known at all nodes on the local Chebyshev grid.
Now perform spectral differentiation to evaluate approximations to $\partial u/\partial x_{2}$
for the Chebyshev nodes on the two horizontal sides, and $\partial u/\partial x_{1}$ for the
Chebyshev nodes on the two vertical sides.

\lsp

\noindent
\textbf{\textit{Step 4 --- re-tabulation from the Chebyshev nodes back to Gaussian nodes:}}
After Step 3, the boundary fluxes on $\partial \Omega_{\tau}$ are specified by four polynomials
of degree $\ngauss-1$ (specified via tabulation on the relevant Chebyshev nodes). Now simply
evaluate these polynomials at the Gaussian nodes on each side to obtain the vector $\vv(I_{\rm e}^{\tau})$.

\lsp

Putting everything together, we find that the matrix $\mtx{T}^{\tau}$ is given as a product of
four matrices
$$
\begin{array}{ccccccccc}
\mtx{T}^{\tau} &=& \mtx{L}_{4} &\circ& \mtx{L}_{3} &\circ& \mtx{L}_{2} &\circ& \mtx{L}_{1} \\
4\ngauss \times 4\ngauss
&&
4\ngauss \times 4\ngauss
&&
4\ngauss \times \ngauss^{2}
&&
\ngauss^{2} \times 4(\ngauss-1)
&&
4(\ngauss-1) \times 4\ngauss
\end{array}
$$
where $\mtx{L}_{i}$ is the linear transform corresponding to ``Step $i$'' above. Observe that
many of these transforms are far from dense, for instance, $\mtx{L}_{1}$ and $\mtx{L}_{4}$
are $4\times 4$ block matrices with all off-diagonal blocks equal to zero. Exploiting these
structures substantially accelerates the computation.

\begin{remark}
The grid of Chebyshev nodes $\{\pzz_{k}\}_{j=1}^{\ncheb^{2}}$ introduced
in Section \ref{sec:prelim} is only used for the local computation. In the final
solver, there is no need to store potential values at these grid points --- they
are used merely for constructing the matrix $\mtx{T}^{\tau}$.
\end{remark}

\section{Merging two DtN maps}
\label{sec:merge}

Let $\tau$ denote a box in the tree with children $\alpha$ and $\beta$. In this
section, we demonstrate that if the DtN matrices $\mtx{T}^{\alpha}$ and $\mtx{T}^{\beta}$
for the children are known, then the DtN matrix $\mtx{T}^{\tau}$ can be constructed
via a purely local computation which we refer to as a ``merge'' operation.

\begin{figure}
\setlength{\unitlength}{1mm}
\begin{picture}(95,55)
\put(-15,00){\includegraphics[height=55mm]{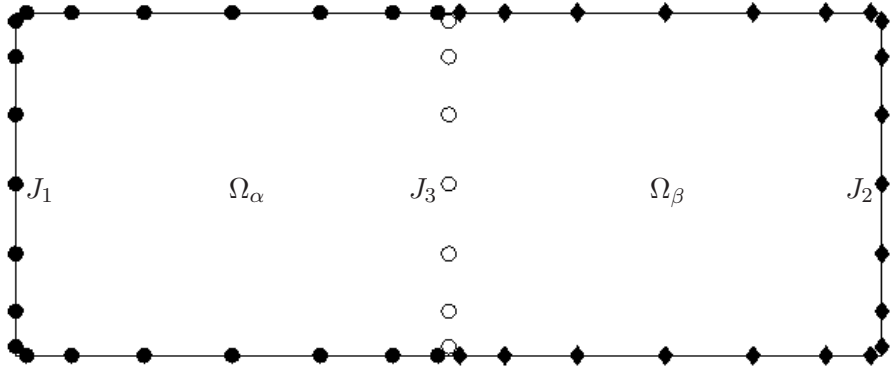}}
\put(18,25){$\Omega_{\alpha}$}
\put(74,25){$\Omega_{\beta}$}
\put(-9,25){$J_{1}$}
\put(100,25){$J_{2}$}
\put(42,25){$J_{3}$}
\end{picture}
\caption{Notation for the merge operation described in Section \ref{sec:merge}.
The rectangular domain $\Omega$ is formed by two squares $\Omega_{\alpha}$ and $\Omega_{\beta}$.
The sets $J_{1}$ and $J_{2}$ form the exterior nodes (black), while
$J_{3}$ consists of the interior nodes (white).}
\label{fig:siblings_notation}
\end{figure}

We start by introducing some notation:
Let $\Omega_{\tau}$ denote a box with children $\Omega_{\alpha}$ and
$\Omega_{\beta}$. For concreteness, let us assume that $\Omega_{\alpha}$ and
$\Omega_{\beta}$ share a vertical edge as shown in Figure \ref{fig:siblings_notation},
so that
$$
\Omega_{\tau} = \Omega_{\alpha} \cup \Omega_{\beta}.
$$
We partition the points on $\partial\Omega_{\alpha}$ and $\partial\Omega_{\beta}$ into three sets:
\begin{tabbing}
\mbox{}\hspace{5mm}\= $J_{1}$ \hspace{4mm} \=
Boundary nodes of $\Omega_{\alpha}$ that are not boundary nodes of $\Omega_{\beta}$.\\
\> $J_{2}$ \> Boundary nodes of $\Omega_{\beta}$ that are not boundary nodes of $\Omega_{\alpha}$.\\
\> $J_{3}$ \> Boundary nodes of both $\Omega_{\alpha}$ and $\Omega_{\beta}$ that are \textit{not} boundary nodes of the union box
$\Omega_{\tau}$.
\end{tabbing}
Figure \ref{fig:siblings_notation} illustrates the definitions of the $J_{k}$'s.
Let $u$ denote a solution to (\ref{eq:basic}), with tabulated potential values $\uu$
and boundary fluxes $\vv$, as described in Section \ref{sec:discretization}. Set
\begin{equation}
\label{eq:defuiue}
\uu_{\rm i} = \uu_{3},
\qquad\mbox{and}\qquad
\uu_{\rm e} = \left[\begin{array}{c} \uu_{1} \\ \uu_{2} \end{array}\right].
\end{equation}
Recall that $\TT^{\alpha}$ and $\TT^{\beta}$ denote the
operators that map values of the potential $u$ on the boundary to values of
$\partial_{n}u$ on the boundaries
of the boxes $\Omega_{\alpha}$ and $\Omega_{\beta}$,
as described in Section \ref{sec:leaf}.
The operators can be partitioned according to the numbering of nodes
in Figure \ref{fig:siblings_notation}, resulting in the equations
\begin{equation}
\label{eq:bittersweet1}
\left[\begin{array}{c}
\vv_{1}\\ \vv_{3}
\end{array}\right] =
\left[\begin{array}{ccc}
\mtx{T}_{1,1}^{\alpha} & \mtx{T}_{1,3}^{\alpha} \\
\mtx{T}_{3,1}^{\alpha} & \mtx{T}_{3,3}^{\alpha}
\end{array}\right]\,
\left[\begin{array}{c}
\uu_{1}\\ \uu_{3}
\end{array}\right],
\qquad {\rm and} \qquad
\left[\begin{array}{c}
\vv_{2}\\ \vv_{3}
\end{array}\right] =
\left[\begin{array}{ccc}
\mtx{T}_{2,2}^{\beta} & \mtx{T}_{2,3}^{\beta} \\
\mtx{T}_{3,2}^{\beta} & \mtx{T}_{3,3}^{\beta}
\end{array}\right]\,
\left[\begin{array}{c}
\uu_{2}\\ \uu_{3}
\end{array}\right].
\end{equation}

Our objective is now to construct a solution operator $\mtx{S}^{\tau}$
and a DtN matrix $\mtx{T}^{\tau}$ such that
\begin{align}
\label{eq:desiredS}
\uu_{3} =&\ \mtx{S}^{\tau}\,
\left[\begin{array}{c}\uu_{1} \\ \uu_{2}\end{array}\right]\\
\label{eq:desiredT}
\left[\begin{array}{c}\vv_{1} \\ \vv_{2}\end{array}\right]
=&\
\mtx{T}^{\tau}\,
\left[\begin{array}{c}\uu_{1} \\ \uu_{2}\end{array}\right].
\end{align}
To this end, we write (\ref{eq:bittersweet1}) as a single equation:
\begin{equation}
\label{eq:premerge}
\left[\begin{array}{cc|c}
\TT^{\alpha}_{1,3} & \mtx{0} & \TT^{\alpha}_{3,3} \\
\mtx{0} & \TT^{\beta }_{2,3} & \TT^{\beta }_{3,3} \\ \hline
\TT^{\alpha}_{1,3} & -\TT^{\beta}_{2,3} & \TT^{\alpha}_{3,3} - \TT^{\beta}_{3,3}
\end{array}\right]\,
\left[\begin{array}{c}
\uu_{1} \\ \uu_{2} \\ \hline \uu_{3}
\end{array}\right]
=
\left[\begin{array}{c}
\vv_{1} \\ \vv_{2} \\ \hline \vct{0}
\end{array}\right],
\end{equation}
The last equation directly tells us that (\ref{eq:desiredS}) holds with
\begin{equation}
\label{eq:U_parent}
\mtx{S}^{\tau} =
\bigl(\TT^{\alpha}_{3,3} - \TT^{\beta}_{3,3}\bigr)^{-1}
\bigl[-\TT^{\alpha}_{3,1}\ \big|\ \TT^{\beta}_{3,2}].
\end{equation}
By eliminating $\uu_{3}$ from (\ref{eq:premerge}) by forming a Schur complement, we also find that
(\ref{eq:desiredT}) holds with
\begin{equation}
\label{eq:merge}
\mtx{T}^{\tau} =
\left[\begin{array}{ccc}
\TT_{1,1}^{\alpha} & \mtx{0} \\
\mtx{0} & \TT_{2,2}^{\beta }
\end{array}\right] +
\left[\begin{array}{c}
\TT_{1,3}^{\alpha} \\
\TT_{2,3}^{\beta}
\end{array}\right]\,
\bigl(\TT^{\alpha}_{3,3} - \TT^{\beta}_{3,3}\bigr)^{-1}
\bigl[-\TT^{\alpha}_{3,1}\ \big|\ \TT^{\beta}_{3,2}\bigr].
\end{equation}


\section{The full hierarchical scheme}
\label{sec:hierarchy}

At this point, we know how to construct the DtN operator for a leaf (Section
\ref{sec:leaf}), and how to merge two such operators of neighboring
patches to form the DtN operator of their union (Section
\ref{sec:merge}). We are ready to describe the full hierarchical
scheme for solving the Dirichlet problem (\ref{eq:basic}). This scheme
takes the Dirichlet boundary data $f$, and constructs an approximation
to the solution $u$. The output is a vector $\uu$ that tabulates
approximations to $u$ at the Gaussian nodes $\{\pxx_{k}\}_{k=1}^{N}$
on all interior edges that were defined in Section \ref{sec:discretization}.
To find $u$ at an arbitrary set of target points in $\Omega$, a post-processing
step described in Section \ref{sec:postproc} can be used.

\subsection{The algorithm}
\label{sec:thealgorithm}
Partition the domain into a hierarchical tree as described in Section \ref{sec:discretization}.
Then execute a ``build stage'' in which we construct for each box $\tau$ the following two matrices:

\vspace{2mm}

\begin{itemize}
\item[$\mtx{S}^{\tau}$]
For a parent box $\tau$, $\mtx{S}^{\tau}$ is a solution operator that maps values
of $u$ on $\partial \Omega_{\tau}$ to values of $u$ at the interior nodes. In other
words, $\uu(I_{\rm i}^{\tau}) = \mtx{S}^{\tau}\,\uu(I_{\rm e}^{\tau})$.
(For a leaf $\tau$, $\mtx{S}^{\tau}$ is not defined.)

\vspace{2mm}

\item[$\TT^{\tau}$] The matrix that maps $\uu(I_{\rm e}^{\tau})$ (tabulating values of $u$ on
$\partial \Omega_{\tau}$) to $\vv(I_{\rm e}^{\tau})$ (tabulating values of $du/dn$).
In other words, $\vv(I_{\rm e}^{\tau}) = \TT^{\tau}\,\uu(I_{\rm e}^{\tau})$.
\end{itemize}

\vspace{2mm}

\noindent
(Recall that the index vector $I_{\rm e}^{\tau}$ and $I_{\rm i}^{\tau}$
were defined in Section \ref{sec:discretization}.)
The build stage consists of a single sweep over all nodes in the tree.
Any bottom-up ordering in which any parent box is processed after its
children can be used. For each leaf box $\tau$, an approximation
to the local DtN map $\mtx{T}^{\tau}$ is constructed using the procedure
described in Section \ref{sec:leaf}. For a parent box
$\tau$ with children $\sigma_{1}$ and $\sigma_{2}$, the matrices
$\mtx{S}^{\tau}$ and $\TT^{\tau}$ are formed from the DtN operators $\TT^{\sigma_{1}}$
and $\TT^{\sigma_{2}}$ via the process described in Section \ref{sec:merge}.
Algorithm 1 summarizes the build stage.


Once all the matrices $\{\mtx{S}^{\tau}\}_{\tau}$ have been formed, a vector $\uu$ holding approximations
to the solution $u$ of (\ref{eq:basic}) can be constructed for all discretization points by
starting at the root box $\Omega$ and moving down the tree toward the leaf boxes.
The values of $\uu$ for the points on boundary of $\Omega$ can be obtained by tabulating the
boundary function $f$.  When any box $\tau$ is processed, the value of $\uu$ is known for all nodes
on its boundary (i.e.~those listed in $I_{\rm e}^{\tau}$). The matrix $\mtx{S}^{\tau}$ directly
maps these values to the values of $\uu$ on the nodes in the interior of $\tau$ (i.e.~those
listed in $I_{\rm i}^{\tau}$). When all nodes have been processed, approximations to $u$ have
constructed for all tabulation nodes on interior edges.
Algorithm 2 summarizes the solve stage.


\begin{figure}
\fbox{
\begin{minipage}{140mm}
\begin{center}
\textsc{Algorithm 1} (build solution operators)
\end{center}

This algorithm builds the global Dirichlet-to-Neumann operator for (\ref{eq:basic}).\\
It also builds all matrices $\mtx{S}^{\tau}$ required for constructing $u$ at any interior point.\\
It is assumed that if node $\tau$ is a parent of node $\sigma$, then $\tau < \sigma$.

\rule{\textwidth}{0.5pt}

\begin{tabbing}
\mbox{}\hspace{7mm} \= \mbox{}\hspace{6mm} \= \mbox{}\hspace{6mm} \= \mbox{}\hspace{6mm} \= \mbox{}\hspace{6mm} \= \kill
(1)\> \textbf{for} $\tau = N_{\rm boxes},\,N_{\rm boxes}-1,\,N_{\rm boxes}-2,\,\dots,\,1$\\
(2)\> \> \textbf{if} ($\tau$ is a leaf)\\
(3)\> \> \> Construct $\mtx{S}^{\tau}$ via the process described in Section \ref{sec:leaf}.\\
(4)\> \> \textbf{else}\\
(5)\> \> \> Let $\sigma_{1}$ and $\sigma_{2}$ be the children of $\tau$.\\
(6)\> \> \> Split $I_{\rm e}^{\sigma_{1}}$ and $I_{\rm e}^{\sigma_{2}}$ into vectors $I_{1}$, $I_{2}$, and $I_{3}$ as shown in Figure \ref{fig:siblings_notation}.\\
(7)\> \> \> $\mtx{S}^{\tau} = \bigl(\TT^{\sigma_{1}}_{3,3} - \TT^{\sigma_{2}}_{3,3}\bigr)^{-1}
                           \bigl[-\TT^{\sigma_{1}}_{3,1}\  \big|\
                                  \TT^{\sigma_{2}}_{3,2}\bigr]$\\
(8)\> \> \> $\TT^{\tau} = \left[\begin{array}{ccc}
                          \mtx{T}_{1,1}^{\sigma_{1}} & \mtx{0}\\
                          \mtx{0} & \mtx{T}_{2,2}^{\sigma_{2} }
                          \end{array}\right] +
                    \left[\begin{array}{c}
                          \TT_{1,3}^{\sigma_{1}} \\
                          \TT_{2,3}^{\sigma_{2}}
                          \end{array}\right]\,\mtx{S}^{\tau}$.\\
(9)\> \> \> Delete $\TT^{\sigma_{1}}$ and $\TT^{\sigma_{1}}$.\\
(10)\> \> \textbf{end if}\\
(11)\> \textbf{end for}
\end{tabbing}
\end{minipage}}
\end{figure}

\begin{figure}
\fbox{
\begin{minipage}{140mm}
\begin{center}
\textsc{Algorithm 2} (solve BVP once solution operator has been built)
\end{center}

This program constructs an approximation $\uu$ to the solution $u$ of (\ref{eq:basic}).\\
It assumes that all matrices $\mtx{S}^{\tau}$ have already been constructed in a pre-computation.\\
It is assumed that if node $\tau$ is a parent of node $\sigma$, then $\tau < \sigma$.

\rule{\textwidth}{0.5pt}

\begin{tabbing}
\mbox{}\hspace{7mm} \= \mbox{}\hspace{6mm} \= \mbox{}\hspace{6mm} \= \mbox{}\hspace{6mm} \= \mbox{}\hspace{6mm} \= \kill
(1)\> $\uu(k) = f(\pxx_{k})$ for all $k \in I_{\rm e}^{1}$.\\
(2)\> \textbf{for} $\tau = 1,\,2,\,3,\,\dots,\,N_{\rm boxes}$\\
(3)\> \> $\uu(I_{\rm i}^{\tau}) = \mtx{S}^{\tau}\,\uu(I_{\rm e}^{\tau})$.\\
(4)\> \textbf{end for}
\end{tabbing}

\vspace{1mm}

\textit{Remark: This algorithm outputs the solution on the Gaussian nodes on box boundaries.
To get the solution at other points, use the method described in Section \ref{sec:postproc}.}

\end{minipage}}
\end{figure}

\begin{remark}
The merge stage is exact when performed in exact arithmetic. The only
approximation involved is the approximation of the solution $u$ on a leaf by its
interpolating polynomial.
\end{remark}

\begin{remark}
\label{remark:refine}
To keep the presentation simple, we consider in this paper only the case
of a uniform computational grid. Such grids are obviously not well suited
to situations where the regularity of the solution changes across the domain.
The method described can \textit{in principle} be modified to handle locally refined
grids quite easily. A complication is that the tabulation nodes for
two touching boxes will typically not coincide, which requires the introduction of
specialized interpolation operators. Efficient refinement strategies also
require the development of error indicators that identify the regions where
the grid need to be refined. This is work in progress, and will be reported
at a later date. We observe that our introduction of Gaussian nodes on the
internal boundaries (as opposed to the Chebyshev nodes used in
\cite{2012_martinsson_spectralcomposite_jcp}) makes re-interpolation much
easier.
\end{remark}

\subsection{Asymptotic complexity}
\label{sec:complexity}
In this section, we determine the asymptotic complexity of the direct solver.
Let $N_{\rm leaf} = 4\ngauss$ denote the number of Gaussian nodes on the boundary
of a leaf box, and let $q^2$ denote the number of Chebychev nodes used in the leaf computation.
Let $L$ denote the number of levels in the binary tree.  This means there are
$4^L$ boxes.  Thus the total number of discretization nodes $N$ is approximately $4^L\ngauss = \frac{(2^L\ngauss)^2}{\ngauss}$.
(To be exact, $N = 2^{2L+1}\ngauss+2^{L+1}\ngauss$.)

The cost to process one leaf is approximately $O(q^6)$.  Since there are $\frac{N}{q^2}$ leaf boxes, the total
cost of pre-computing approximate DtN operators for all the bottom level is
$\frac{N}{q^2}\times q^6\sim N\ngauss^4$.

Next, consider the cost of constructing the DtN map on level $\ell$ via
the merge operation described in Section \ref{sec:merge}.  For each box
on the level $\ell$, the operators $\TT^\tau$ and $\mtx{S}^\tau$ are constructed
via (\ref{eq:U_parent}) and (\ref{eq:U_parent}).  These operations involve
matrices of size roughly $2^{-\ell}N^{0.5}\times  2^{-\ell}N^{0.5}$.  Since there are
$4^\ell$ boxes per level.  The cost on level $\ell$ of the merge is
$$4^\ell\times \left(2^{-\ell}N^{0.5}\right)^3 \sim 2^{-\ell}N^{1.5}.$$

The total cost for all the merge procedures has complexity
$$\sum_{\ell=1}^L 2^{-\ell}N^{1.5} \sim N^{1.5}.$$

Finally, consider the cost of the downwards sweep which solves for the
interior unknowns.  For any non-leaf box $\tau$ on level $\ell$, the size of $\mtx{S}^\tau$ is
$2^l\ngauss \times  2^l(6\ngauss)$ which is approximately $\sim 2^{-\ell}N^{0.5} \times 2^{-\ell}N^{0.5}$.
 Thus the cost of applying $\mtx{S}^\tau$ is roughly $(2^{-\ell}N^{0.5})^2 = 2^{-2\ell}N$.
So the total cost of the solve step has complexity
$$ \sum_{l=0}^{L-1}2^{2\ell}2^{-2\ell}N \sim N\log N.$$

In Section \ref{sec:accel}, we explain how to exploit structure in
the matrices $\TT$ and $\mtx{S}$ to improve the computational
cost of both the precomputation and the solve steps.

\subsection{Post-processing}
\label{sec:postproc}

The direct solver in Algorithm 1 constructs approximations
to the solution $u$ of (\ref{eq:basic}) at tabulation nodes at all interior edges.
Once these are available, it is easy to construct an approximation to $u$ at
an arbitrary point. To illustrate the process, suppose that we seek an approximation
to $u(\pyy)$, where $\pyy$ is a point located in a leaf $\tau$. We have values of $u$
tabulated at Gaussian nodes on $\partial \Omega_{\tau}$. These can easily be re-interpolated
to the Chebyshev nodes on $\partial \Omega_{\tau}$. Then $u$ can be reconstructed at
the interior Chebyshev nodes via the formula (\ref{eq:leaf_solve}); observe that
the local solution operator $-\mtx{A}_{\rm i,i}^{-1}\mtx{A}_{\rm i,e}$ was built when
the leaf was originally processed and can be simply retrieved from memory (assuming
enough memory is available). Once $u$ is tabulated at the Chebyshev grid on $\Omega_{\tau}$, it is
trivial to interpolate it to $\pyy$ or any other point.

\section{Compressible matrices}
\label{sec:HBS}

The cost of the direct solver given as Algorithm 1
is dominated by the work done at the very top levels;
the matrix operations on lines (7) and (8) involve dense
matrices of size $O(N^{0.5}) \times O(N^{0.5})$ where $N$
is the total number of discretization nodes, resulting in
$O(N^{1.5})$ overall cost. It turns out that these dense
matrices have internal structure that can be exploited to
greatly accelerate the matrix algebra. Specifically,
the off-diagonal blocks of these matrices are to high
precision rank deficient, and the matrices can be represented efficiently
using a hierarchical ``data-sparse'' format known as
\textit{Hierarchically Block Separable (HBS)}
(and sometimes \textit{Hierarchically Semi-Separable (HSS)} matrices \cite{2007_shiv_sheng,2004_gu_divide}).
This section briefly describes the HBS property,
for details see \cite{m2011_1D_survey}.

\subsection{Block separable}
Let $\mtx{H}$ be an $mp\times mp$ matrix that is blocked into $p\times p$ blocks,
each of size $m\times m$.
We say that $\mtx{H}$ is ``block separable'' with ``block-rank'' $k$
if for $\tau = 1,\,2,\,\dots,\,p$, there exist $m\times k$
matrices $\mtx{U}_{\tau}$ and $\mtx{V}_{\tau}$ such that each off-diagonal
block $\mtx{H}_{\sigma,\tau}$ of $\mtx{H}$ admits the factorization
\begin{equation}
\label{eq:yy1}
\begin{array}{cccccccc}
\mtx{H}_{\sigma,\tau}  & = & \mtx{U}_{\sigma}   & \tilde{\mtx{H}}_{\sigma,\tau}  & \mtx{V}_{\tau}^{*}, &
\quad \sigma,\tau \in \{1,\,2,\,\dots,\,p\},\quad \sigma \neq \tau.\\
m\times m &   & m\times k & k \times k & k\times m
\end{array}
\end{equation}
Observe that the columns of $\mtx{U}_{\sigma}$ must form a basis for
the columns of all off-diagonal blocks in row $\sigma$, and
analogously, the columns of $\mtx{V}_{\tau}$ must form a basis for the
rows in all the off-diagonal blocks in column $\tau$. When (\ref{eq:yy1})
holds, the matrix $\mtx{H}$ admits a block factorization
\begin{equation}
\label{eq:yy2}
\begin{array}{cccccccccc}
 \mtx{H} &  =& \mtx{U}&\tilde{\mtx{H}}& \mtx{V}^{*} & +&  \mtx{D},\\
mp\times mp &   & mp\times kp & kp \times kp & kp\times mp && mp \times mp\\
\end{array}
\end{equation}
where
$$
\mtx{U} = \mbox{diag}(\mtx{U}_{1},\,\mtx{U}_{2},\,\dots,\,\mtx{U}_{p}),\quad
\mtx{V} = \mbox{diag}(\mtx{V}_{1},\,\mtx{V}_{2},\,\dots,\,\mtx{V}_{p}),\quad
\mtx{D} = \mbox{diag}(\mtx{D}_{1},\,\mtx{D}_{2},\,\dots,\,\mtx{D}_{p}),
$$
and
$$\tilde{\mtx{H}} = \left[\begin{array}{cccc}
\mtx{0} & \tilde{\mtx{H}}_{12} & \tilde{\mtx{H}}_{13} & \cdots \\
\tilde{\mtx{H}}_{21} & \mtx{0} & \tilde{\mtx{H}}_{23} & \cdots \\
\tilde{\mtx{H}}_{31} & \tilde{\mtx{H}}_{32} & \mtx{0} & \cdots \\
\vdots & \vdots & \vdots
\end{array}\right].
$$

\subsection{Hierarchically Block-Separable}
Informally speaking, a matrix $\mtx{H}$ is \textit{Heirarchically Block-Separable} (HBS),
if it is amenable to a \textit{telescoping} block factorization.  In other words,
in addition to the matrix $\mtx{H}$ being block separable, so is $\tilde{\mtx{H}}$
once it has been reblocked to form a matrix with $p/2 \times p/2$ blocks.
Likewise, the middle matrix from the block separable factorization
of $\tilde{\mtx{H}}$ will be block separable, etc.

In this section, we describe properties and the factored representation of HBS matrices.
Details on constructing the factorization are provided in \cite{m2011_1D_survey}.

\subsection{A binary tree structure}
\label{sec:tree}

\begin{figure}
\setlength{\unitlength}{1mm}
\begin{picture}(169,41)
\put(20, 0){\includegraphics[height=41mm]{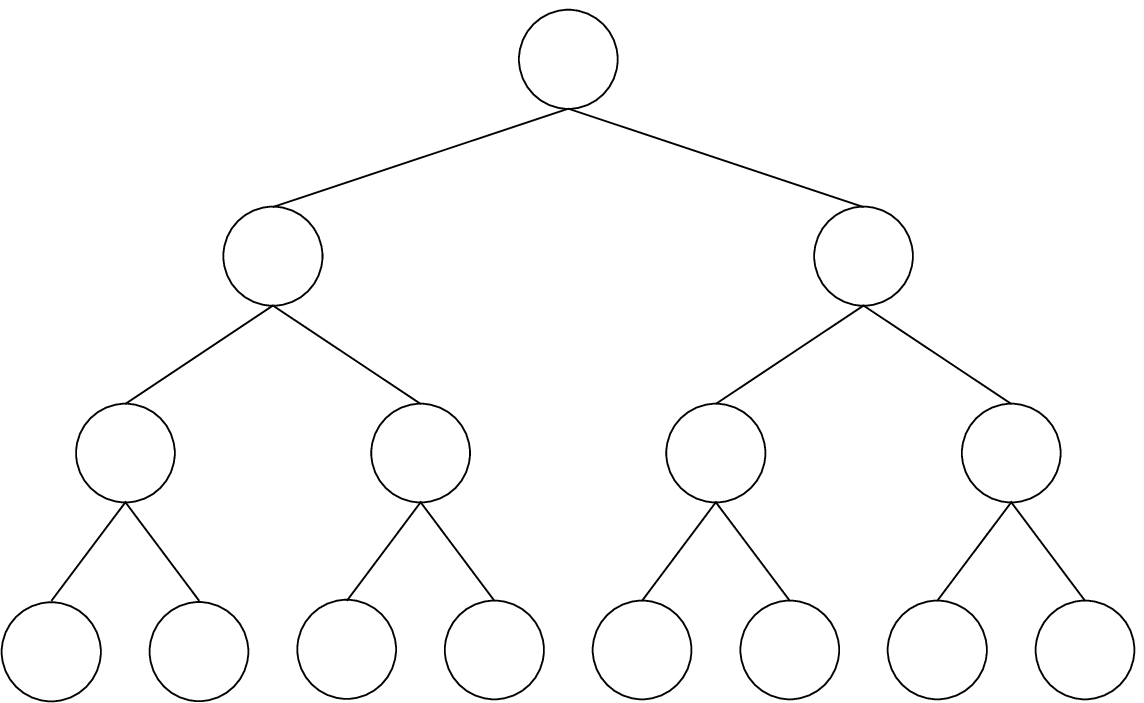}}
\put( 0,36){Level $0$}
\put( 0,25){Level $1$}
\put( 0,13){Level $2$}
\put( 0, 2){Level $3$}
\put(90,36){\footnotesize$I_{1} = [1,\,2,\,\dots,\,400]$}
\put(90,25){\footnotesize$I_{2} = [1,\,2,\,\dots,\,200]$, $I_{3} = [201,\,202,\,\dots,\,400]$}
\put(90,13){\footnotesize$I_{4} = [1,\,2,\,\dots,\,100]$, $I_{5} = [101,\,102,\,\dots,\,200]$, \dots}
\put(90, 2){\footnotesize$I_{8} = [1,\,2,\,\dots,\,50]$, $I_{9} = [51,\,52,\,\dots,\,100]$, \dots}
\put(52.5,36.5){$1$}
\put(35,25){$2$}
\put(70,25){$3$}
\put(26,13){$4$}
\put(44,13){$5$}
\put(61,13){$6$}
\put(78.5,13){$7$}
\put(22, 2){\small$8$}
\put(31, 2){\small$9$}
\put(38.5, 2){\small$10$}
\put(47, 2){\small$11$}
\put(56, 2){\small$12$}
\put(64.5, 2){\small$13$}
\put(73, 2){\small$14$}
\put(82, 2){\small$15$}
\end{picture}
\caption{Numbering of nodes in a fully populated binary tree with $L=3$ levels.
The root is the original index vector $I = I_{1} = [1,\,2,\,\dots,\,400]$.}
\label{fig:tree}
\end{figure}

The HBS representation
of an $M\times M$ matrix $\mtx{H}$ is
based on a partition of the index vector $I = [1,\,2,\,\dots,\,M]$
into a binary tree structure.
We let $I$ form the root of the tree, and give it the index $1$,
$I_{1} = I$. We next split the root into two roughly equi-sized
vectors $I_{2}$ and $I_{3}$ so that $I_{1} = I_{2} \cup I_{3}$.
The full tree is then formed by continuing to subdivide any interval
that holds more than some preset fixed number $m$ of indices.
We use the integers $\ell = 0,\,1,\,\dots,\,L$ to label the different
levels, with $0$ denoting the coarsest level.
A \textit{leaf} is a node corresponding to a vector that never got split.
For a non-leaf node $\tau$, its \textit{children} are the two boxes
$\sigma_{1}$ and $\sigma_{2}$ such that $I_{\tau} = I_{\sigma_{1}} \cup I_{\sigma_{2}}$,
and $\tau$ is then the \textit{parent} of $\sigma_{1}$ and $\sigma_{2}$.
Two boxes with the same parent are called \textit{siblings}. These
definitions are illustrated in Figure \ref{fig:tree}.

\subsection{Definition of the HBS property}
\label{sec:rigor}
We now define what it means for an $M\times M$ matrix $\mtx{H}$
to be \textit{hierarchically block separable} with respect to a given binary tree $\mathcal{T}$
that partitions the index vector $J = [1,\,2,\,\dots,\,M]$. For simplicity, we suppose that for every leaf node $\tau$ the index vector
$I_{\tau}$ holds precisely $m$ points, so that $M = m\,2^{L}$. Then $\mtx{H}$ is HBS with block
rank $k$ if the following two conditions hold:

\lsp

\noindent
\textit{(1) Assumption on ranks of off-diagonal blocks at the finest level:}
For any two distinct leaf nodes $\tau$ and $\tau'$, define the $m\times m$ matrix
\begin{equation}
\label{eq:def1}
\mtx{H}_{\tau,\tau'} = \mtx{H}(I_{\tau},I_{\tau'}).
\end{equation}
Then there must exist matrices $\mtx{U}_{\tau}$, $\mtx{V}_{\tau'}$, and $\tilde{\mtx{H}}_{\tau,\tau'}$ such that
\begin{equation}
\label{eq:def2}
\begin{array}{cccccccccc}
\mtx{H}_{\tau,\tau'} & = & \mtx{U}_{\tau} & \tilde{\mtx{H}}_{\tau,\tau'} & \mtx{V}_{\tau'}^{*}.  \\
m\times m &   & m\times k & k \times k & k\times m
\end{array}
\end{equation}

\lsp

\noindent
\textit{(2) Assumption on ranks of off-diagonal blocks on level $\ell = L-1,\,L-2,\,\dots,\,1$:}
The rank assumption at level $\ell$ is defined in terms of the blocks constructed
on the next finer level $\ell+1$: For any distinct nodes $\tau$ and
$\tau'$ on level $\ell$ with  children $\sigma_{1},\sigma_{2}$ and $\sigma_{1}',\sigma_{2}'$,
respectively, define
\begin{equation}
\label{eq:def3}
\mtx{H}_{\tau,\tau'} = \left[\begin{array}{cc}
\tilde{\mtx{H}}_{\sigma_{1},\sigma_{1}'} & \tilde{\mtx{H}}_{\sigma_{1},\sigma_{2}'} \\
\tilde{\mtx{H}}_{\sigma_{2},\sigma_{1}'} & \tilde{\mtx{H}}_{\sigma_{2},\sigma_{2}'}
\end{array}\right].
\end{equation}
Then there must exist matrices $\mtx{U}_{\tau}$, $\mtx{V}_{\tau'}$, and $\tilde{\mtx{H}}_{\tau,\tau'}$ such that
\begin{equation}
\label{eq:def4}
\begin{array}{cccccccccc}
\mtx{H}_{\tau,\tau'} & = & \mtx{U}_{\tau} & \tilde{\mtx{H}}_{\tau,\tau'} & \mtx{V}_{\tau'}^{*}.  \\
2k\times 2k &   & 2k\times k & k \times k & k\times 2k
\end{array}
\end{equation}

\lsp

 An HBS matrix is
now fully described if the basis matrices $\mtx{U}_{\tau}$ and $\mtx{V}_{\tau}$ are provided
for each node $\tau$, and in addition, we are for each leaf $\tau$ given the $m\times m$ matrix
\begin{equation}
\label{eq:def5}
\mtx{D}_{\tau} = \mtx{H}(I_{\tau},I_{\tau}),
\end{equation}
and for each parent node $\tau$ with children $\sigma_{1}$ and $\sigma_{2}$ we are given
the $2k\times 2k$ matrix
\begin{equation}
\label{eq:def6}
\mtx{B}_{\tau}       = \left[\begin{array}{cc}
0 & \tilde{\mtx{H}}_{\sigma_{1},\sigma_{2}} \\
\tilde{\mtx{H}}_{\sigma_{2},\sigma_{1}} & 0
\end{array}\right].
\end{equation}
Observe in particular that the matrices $\tilde{\mtx{H}}_{\sigma_{1},\sigma_{2}}$ are only required when
$\{\sigma_{1},\sigma_{2}\}$
forms a sibling pair. Figure \ref{fig:summary_of_factors} summarizes the required matrices.

\begin{figure}
\small
\begin{tabular}{|l|l|l|l|} \hline
                & Name:            & Size:       & Function: \\ \hline
For each leaf   & $\mtx{D}_{\tau}$ & $ m\times  m$ & The diagonal block $\mtx{H}(I_{\tau},I_{\tau})$. \\
node $\tau$:    & $\mtx{U}_{\tau}$ & $ m\times  k$ & Basis for the columns in the blocks in row $\tau$. \\
                & $\mtx{V}_{\tau}$ & $ m\times  k$ & Basis for the rows in the blocks in column $\tau$. \\ \hline
For each parent & $\mtx{B}_{\tau}$ & $2k\times 2k$ & Interactions between the children of $\tau$. \\
node $\tau$:    & $\mtx{U}_{\tau}$ & $2k\times  k$ & Basis for the columns in the (reduced) blocks in row $\tau$. \\
                & $\mtx{V}_{\tau}$ & $2k\times  k$ & Basis for the rows in the (reduced) blocks in column $\tau$. \\ \hline
\end{tabular}
\caption{An HBS matrix $\mtx{H}$ associated with a tree $\mathcal{T}$ is fully specified if the factors
listed above are provided.}
\label{fig:summary_of_factors}
\end{figure}

\subsection{Telescoping factorization}
Given the matrices defined in the previous section,  we define the following block diagonal factors:
\begin{align}
\label{eq:def_uD}
\underline{\mtx{D}}^{(\ell)} &= \mbox{diag}(\mtx{D}_{\tau}\,\colon\, \tau\mbox{ is a box on level }\ell),\qquad \ell = 0,\,1,\,\dots,\,L,\\
\underline{\mtx{U}}^{(\ell)} &= \mbox{diag}(\mtx{U}_{\tau}\,\colon\, \tau\mbox{ is a box on level }\ell),\qquad \ell = 1,\,2,\,\dots,\,L,\\
\underline{\mtx{V}}^{(\ell)} &= \mbox{diag}(\mtx{V}_{\tau}\,\colon\, \tau\mbox{ is a box on level }\ell),\qquad \ell = 1,\,2,\,\dots,\,L,\\
\underline{\mtx{B}}^{(\ell)} &= \mbox{diag}(\mtx{B}_{\tau}\,\colon\, \tau\mbox{ is a box on level }\ell),\qquad \ell = 0,\,1,\,\dots,\,L-1,.
\end{align}
 Furthermore, we let $\tilde{\mtx{H}}^{(\ell)}$ denote the block matrix whose diagonal blocks are zero,
and whose off-diagonal blocks are the blocks $\tilde{\mtx{H}}_{\tau,\tau'}$ for all distinct $\tau,\tau'$
on level $\ell$.
With these definitions,
\begin{equation}
\label{eq:tele1}
\begin{array}{cccccccccccc}
\mtx{H} & = & \underline{\mtx{U}}^{(L)} & \tilde{\mtx{H}}^{(L)} & (\underline{\mtx{V}}^{(L)})^{*} & + & \underline{\mtx{D}}^{(L)};\\
m\,2^{L} \times n\,2^{L} &&
m\,2^{L} \times k\,2^{L} &
k\,2^{L} \times k\,2^{L} &
k\,2^{L} \times m\,2^{L} &&
m\,2^{L} \times m\,2^{L}
\end{array}
\end{equation}
for $\ell = L-1,\,L-2,\,\dots,\,1$ we have
\begin{equation}
\label{eq:tele2}
\begin{array}{cccccccccccc}
\tilde{\mtx{H}}^{(\ell+1)} &=& \underline{\mtx{U}}^{(\ell)} & \tilde{\mtx{H}}^{(\ell)} & (\underline{\mtx{V}}^{(\ell)})^{*} & +
& \underline{\mtx{B}}^{(\ell)};\\
k\,2^{\ell+1} \times k\,2^{\ell+1} &&
k\,2^{\ell+1} \times k\,2^{\ell} &
k\,2^{\ell} \times k\,2^{\ell} &
k\,2^{\ell} \times k\,2^{\ell+1} &&
k\,2^{\ell+1} \times k\,2^{\ell+1}
\end{array}
\end{equation}
and finally
\begin{equation}
\label{eq:tele3}
\tilde{\mtx{H}}^{(1)} = \underline{\mtx{B}}^{(0)}.
\end{equation}

\section{Fast arithmetic operations on HBS matrices}
\label{sec:fastHBS}

Arithmetic operations involving dense HBS matrices of size $M\times M$ can often
be executed in $O(M)$ operations. This fast matrix algebra is vital for achieving
linear complexity in our direct solver. This section provides a brief introduction
to the HBS matrix algebra. We describe the operations we need (inversion, addition,
and low-rank update) in some detail for the single level ``block separable'' format.
The generalization to the multi-level ``hierarchically block separable'' format is
briefly described for the case of matrix inversion. A full description of all
algorithms required is given in \cite{Adiss}, which is related to the earlier work
\cite{2010_jianlin_fast_hss}.

Before we start, we recall that a block separable matrix $\mtx{H}$ consisting of
$p\times p$ blocks, each of size $m\times m$, and with ``HBS-rank'' $k<m$, admits
the factorization
\begin{equation}
\label{eq:yy2b}
\begin{array}{cccccccccc}
 \mtx{H} &  =& \mtx{U}&\tilde{\mtx{H}}& \mtx{V}^{*} & +&  \mtx{D}.\\
mp\times mp &   & mp\times kp & kp \times kp & kp\times mp && mp \times mp\\
\end{array}
\end{equation}

\subsection{Inversion of a block separable matrix}
The decomposition (\ref{eq:yy2b}) represents $\mtx{H}$ as a sum of
one term $\mtx{U}\tilde{\mtx{H}}\mtx{V}^{*}$ that is ``low rank,''
and one term $\mtx{D}$ that is easily invertible (since it is block diagonal).
By modifying the classical Woodbury formula for inversion of a matrix
perturbed by the addition of a low-rank term, it can be shown that (see Lemma 3.1 of
\cite{m2011_1D_survey})
\begin{equation}
\label{eq:woodbury}
\mtx{H}^{-1} = \mtx{E}\,(\tilde{\mtx{H}} + \hat{\mtx{D}})^{-1}\,\mtx{F}^{*} + \mtx{G},
\end{equation}
where
\begin{align}
\label{eq:def_muhD}
\hat{\mtx{D}} =&\ \bigl(\mtx{V}^{*}\,\mtx{D}^{-1}\,\mtx{U}\bigr)^{-1},\\
\label{eq:def_muE}
\mtx{E}  =&\ \mtx{D}^{-1}\,\mtx{U}\,\hat{\mtx{D}},\\
\label{eq:def_muF}
\mtx{F}  =&\ (\hat{\mtx{D}}\,\mtx{V}^{*}\,\mtx{D}^{-1})^{*},\\
\label{eq:def_muG}
\mtx{G}  =&\ \mtx{D}^{-1} - \mtx{D}^{-1}\,\mtx{U}\,\hat{\mtx{D}}\,\mtx{V}^{*}\,\mtx{D}^{-1},
\end{align}
assuming the inverses in formulas (\ref{eq:woodbury}) --- (\ref{eq:def_muG}) all exist.
Now observe that the matrices $\hat{\mtx{D}}$, $\mtx{E}$, $\mtx{F}$, and $\mtx{G}$
can all easily be computed since the formulas defining them involve only block-diagonal matrices.
In consequence, (\ref{eq:woodbury}) reduces the task of inverting the big (size $mp\times mp$)
matrix $\mtx{H}$ to the task of inverting the small (size $kp\times kp$) matrix $\tilde{\mtx{H}} + \hat{\mtx{D}}$.

When $\mtx{H}$ is not only ``block separable'', but ``hierarchically block separable'', the
process can be repeated recursively by exploiting that $\tilde{\mtx{H}} + \hat{\mtx{D}}$ is
itself amenable to accelerated inversion, etc. The resulting process is somewhat tricky to
analyze, but leads to very clean codes. To illustrate, we include Algorithm 3
which shows the multi-level $O(M)$ inversion algorithm for an HBS matrix $\mtx{H}$.
The algorithm takes as input the factors $\{\mtx{U}_{\tau},\,\mtx{V}_{\tau},\,\mtx{D}_{\tau},\,\mtx{B}_{\tau}\}_{\tau}$
representing $\mtx{H}$ (cf.~Figure \ref{fig:summary_of_factors}), and outputs an analogous
set of factors $\{\mtx{E}_{\tau},\,\mtx{F}_{\tau},\,\mtx{G}_{\tau}\}_{\tau}$ representing
$\mtx{H}^{-1}$. With these factors, the matrix-vector multiplication $\vct{y} = \mtx{H}^{-1}\vct{x}$
can be executed via the procedure described in Algorithm 4.

\begin{figure}[ht]
\begin{center}
\fbox{
\begin{minipage}{.9\textwidth}\small

\begin{center}
\textsc{Algorithm 3} (inversion of an HBS matrix)
\end{center}

Given factors $\{\mtx{U}_{\tau},\,\mtx{V}_{\tau},\,\mtx{D}_{\tau},\,\mtx{B}_{\tau}\}_{\tau}$
representing an HBS matrix $\mtx{H}$, this algorithm constructs factors
$\{\mtx{E}_{\tau},\,\mtx{F}_{\tau},\,\mtx{G}_{\tau}\}_{\tau}$ representing
$\mtx{H}^{-1}$.

\begin{tabbing}
\hspace{5mm} \= \hspace{5mm} \= \hspace{5mm} \= \kill
\textbf{loop} over all levels, finer to coarser, $\ell = L,\,L-1,\,\dots,1$\\
\> \textbf{loop} over all boxes $\tau$ on level $\ell$,\\
\> \> \textbf{if} $\tau$ is a leaf node\\
\> \> \> $\tilde{\mtx{D}}_{\tau} = \mtx{D}_{\tau}$\\
\> \> \textbf{else}\\
\> \> \> Let $\sigma_{1}$ and $\sigma_{2}$ denote the children of $\tau$.\\
\> \> \> $\tilde{\mtx{D}}_{\tau} = \mtwo{\hat{\mtx{D}}_{\sigma_{1}}}{\mtx{B}_{\sigma_{1},\sigma_{2}}}{\mtx{B}_{\sigma_{2},\sigma_{1}}}{\hat{\mtx{D}}_{\sigma_{2}}}$\\
\> \> \textbf{end if}\\
\> \> $\hat{\mtx{D}}_{\tau} = \bigl(\mtx{V}_{\tau}^{*}\,\tilde{\mtx{D}}_{\tau}^{-1}\,\mtx{U}_{\tau}\bigr)^{-1}$.\\
\> \> $\mtx{E}_{\tau} = \tilde{\mtx{D}}_{\tau}^{-1}\,\mtx{U}_{\tau}\,\hat{\mtx{D}}_{\tau}$.\\
\> \> $\mtx{F}_{\tau}^{*} = \hat{\mtx{D}}_{\tau}\,\mtx{V}_{\tau}^{*}\,\tilde{\mtx{D}}_{\tau}^{-1}$.\\
\> \> $\mtx{G}_{\tau} = \tilde{\mtx{D}}_{\tau}^{-1} - \tilde{\mtx{D}}_{\tau}^{-1}\,\mtx{U}_{\tau}\,\hat{\mtx{D}}_{\tau}\,\mtx{V}_{\tau}^{*}\,\tilde{\mtx{D}}_{\tau}^{-1}$.\\
\> \textbf{end loop}\\
\textbf{end loop}\\
$\mtx{G}_{1} = \mtwo{\hat{\mtx{D}}_{2}}{\mtx{B}_{2,3}}{\mtx{B}_{3,2}}{\hat{\mtx{D}}_{3}}^{-1}$.
\end{tabbing}
\end{minipage}}
\end{center}
\end{figure}

\begin{figure}
\begin{center}
\fbox{
\begin{minipage}{.9\textwidth}\small

\begin{center}
\textsc{Algorithm 4} (application of the inverse of an HBS matrix)
\end{center}

\textit{Given $\vct{x}$, compute $\vct{y} = \mtx{H}^{-1}\,\vct{x}$ using the
factors $\{\mtx{E}_{\tau},\,\mtx{F}_{\tau},\,\mtx{G}_{\tau}\}_{\tau}$
resulting from Algorithm 3.}

\begin{tabbing}
\hspace{5mm} \= \hspace{5mm} \= \hspace{5mm} \= \kill
\textbf{loop} over all leaf boxes $\tau$\\
\> $\hat{\vct{x}}_{\tau} = \mtx{F}_{\tau}^{*}\,\vct{x}(I_{\tau})$.\\
\textbf{end loop}\\[1mm]
\textbf{loop} over all levels, finer to coarser, $\ell = L,\,L-1,\,\dots,1$\\
\> \textbf{loop} over all parent boxes $\tau$ on level $\ell$,\\
\> \> Let $\sigma_{1}$ and $\sigma_{2}$ denote the children of $\tau$.\\
\> \> $\hat{\vct{x}}_{\tau} = \mtx{F}_{\tau}^{*}\,\vtwo{\hat{\vct{x}}_{\sigma_{1}}}{\hat{\vct{x}}_{\sigma_{2}}}$.\\
\> \textbf{end loop}\\
\textbf{end loop}\\[1mm]
$\vtwo{\hat{\vct{y}}_{2}}{\hat{\vct{y}}_{3}} = \mtx{\mtx{G}}_{1}\,\vtwo{\hat{\vct{x}}_{2}}{\hat{\vct{x}}_{3}}$.\\[1mm]
\textbf{loop} over all levels, coarser to finer, $\ell = 1,\,2,\,\dots,\,L-1$\\
\> \textbf{loop} over all parent boxes $\tau$ on level $\ell$\\
\> \> Let $\sigma_{1}$ and $\sigma_{2}$ denote the children of $\tau$.\\
\> \> $\vtwo{\hat{\vct{y}}_{\sigma_{1}}}{\hat{\vct{y}}_{\sigma_{2}}} =
       \mtx{E}_{\tau}\,\hat{\vct{x}}_{\tau} +
       \mtx{G}_{\tau}\,\vtwo{\hat{\vct{x}}_{\sigma_{1}}}{\hat{\vct{x}}_{\sigma_{2}}}$.\\
\> \textbf{end loop}\\
\textbf{end loop}\\[1mm]
\textbf{loop} over all leaf boxes $\tau$\\
\> $\vct{y}(I_{\tau}) = \mtx{E}_{\tau}\,\hat{\vct{q}}_{\tau} + \mtx{G}_{\tau}\,\vct{x}(I_{\tau})$.\\
\textbf{end loop}
\end{tabbing}
\end{minipage}}
\end{center}
\end{figure}


\subsection{Addition of two block separable matrices}
Let $\mtx{H}^A$ and $\mtx{H}^B$ be block separable matrices with factorizations
$$
\mtx{H}^A   = \mtx{U}^A \tilde{\mtx{H}}^A \mtx{V}^{A*}  +  \mtx{D}^A,
\qquad\mbox{and}\qquad
\mtx{H}^B   = \mtx{U}^B \tilde{\mtx{H}}^B \mtx{V}^{B*}  +  \mtx{D}^B.
$$
Then $\mtx{H} = \mtx{H}^{A}+\mtx{H}^{B}$ can be written in block separable form via
\begin{equation}
\label{eq:postadd}
\mtx{H} = \mtx{H}^{A}+\mtx{H}^{B} = \left[\mtx{U}^A\, \mtx{U}^B\right]
\mtwo{\tilde{\mtx{H}}^A}{0}{0}{\tilde{\mtx{H}}^B} \left[\mtx{V}^A\, \mtx{V}^B\right]^*
+\left(\mtx{D}^A+\mtx{D}^B\right).
\end{equation}
To restore (\ref{eq:postadd}) to block separable form, permute the rows and columns
of $\left[\mtx{U}^A\, \mtx{U}^B\right]$ and $\left[\mtx{V}^A\, \mtx{V}^B\right]$ to
attain block diagonal form, then re-orthogonalize the diagonal blocks. This process
in principle results in a matrix $\mtx{H}$ whose HBS-rank is the sum of the HBS-ranks
of $\mtx{H}^{A}$ and $\mtx{H}^{B}$. In practice, this rank increase can be combated
by numerically recompressing the basis matrices, and updating the middle factor as
needed. For details, as well as the extension to a multi-level scheme, see
\cite{2010_jianlin_fast_hss,Adiss}.

\subsection{Addition of a block separable matrix with a low rank matrix}
Let $\mtx{H}^{B}= \mtx{Q}\mtx{R}$ be a $k$-rank matrix where $\mtx{Q}$ and $\mtx{R}^*$ are of size $mp\times k$.  We
would like to add $\mtx{H}^B$ to the block separable matrix $\mtx{H}^A$.
Since we already know how to add two block separable matrices, we choose
to rewrite $\mtx{H}^B$ in block separable form.  Without
loss of generality, assume $\mtx{Q}$ is orthogonal.  Partition $\mtx{Q}$ into $p$ blocks of size $m\times k$.
The blocks make up the matrix $\mtx{U}^B$.  Likewise partition $\mtx{R}$ into $p$ blocks of size $k \times m$.
The block matrix $\mtx{D}^B$ has entries $\mtx{D}_{\tau} = \mtx{Q}_{\tau}\mtx{R}_{\tau}$ for $\tau = 1,\ldots,p$.
  To construct the matrices $\mtx{V}^B$, for each $\tau = 1,\ldots,p$, the matrix $\mtx{R}_{\tau}$ is
factorized into $\tilde{\mtx{R}}_{\tau} \mtx{V}_{\tau}*$ where the matrix $\mtx{V}_{\tau}$ is
orthogonal. The matrices $\tilde{\mtx{R}}_{\tau}$ make up the entries of $\tilde{\mtx{H}}^B$.

\section{Accelerating the direct solver}
\label{sec:accel}

This section describes how the fast matrix algebra described in Sections
\ref{sec:HBS} and \ref{sec:fastHBS} can be used to accelerate the direct
solver of Section \ref{sec:hierarchy} to attain $O(N)$ complexity.
We recall that the $O(N^{1.5})$ cost of Algorithm 1 relates to the execution
of lines (7) and (8) at the top levels, since these involve
dense matrix algebra of matrices of size $O(N^{0.5}) \times O(N^{0.5})$.
The principal claims of this section are:
\begin{itemize}
\item The matrices
$\mtx{T}_{1,3}^{\sigma_{1}}$, $\mtx{T}_{3,1}^{\sigma_{1}}$,
$\mtx{T}_{2,3}^{\sigma_{2}}$, $\mtx{T}_{3,2}^{\sigma_{2}}$
have low numerical rank.
\item The matrices
$\mtx{T}_{1,1}^{\sigma_{1}}$,
$\mtx{T}_{2,2}^{\sigma_{2}}$,
$\mtx{T}_{3,3}^{\sigma_{1}}$,
$\mtx{T}_{3,3}^{\sigma_{2}}$
are HBS matrices of low HBS rank.
\end{itemize}
To be precise, the ranks that we claim are ``low'' scale as
$\log(1/\nu)\times\log(m)$ where $m$ is the number of points along
the boundary of $\Omega_{\tau}$, and $\nu$ is the computational
tolerance. In practice, we found that for problems with non-oscillatory
solutions, the ranks are extremely modest: when $\nu = 10^{-10}$, the
ranks range between $10$ and $80$, even for very large problems.

The cause of the rank deficiencies is that the matrix $\mtx{T}^{\tau}$
is a highly accurate approximation to the Dirichlet-to-Neumann operator
on $\Omega_{\tau}$. This operator is known to have a smooth kernel that
is non-oscillatory whenever the underlying PDE has non-oscillatory solutions.
Since the domain boundary $\partial \Omega_{\tau}$ is one-dimensional, this makes
the expectation that the off-diagonal blocks have low rank very natural,
see \cite{m2011_1D_survey}. It is backed up by extensive numerical
experiments (see Section \ref{sec:numerics}), but we do not at this point
have rigorous proofs to support the claim.

Once it is observed that all matrices in lines (7) and (8) of Algorithm 1
are structured, it becomes obvious how to accelerate the algorithm. For instance,
line (7) is executed in three steps:
(i) Add the HBS matrices $\mtx{T}_{3,3}^{\sigma_{1}}$ and $-\mtx{T}_{3,3}^{\sigma_{2}}$.
(ii) Invert the sum of the HBS matrices.
(iii) Apply the inverse (in HBS form) to one of the low rank factors of
$\bigl[-\TT^{\alpha}_{3,1}\ \big|\ \TT^{\beta}_{3,2}\bigr]$. The result is
an approximation to $\mtx{S}^{\tau}$, represented as a product of two thin
matrices. Executing line (8) is analogous: First form the matrix products
$\mtx{T}_{1,3}^{\sigma_{1}}\,\mtx{S}^{\tau}$ and
$\mtx{T}_{2,3}^{\sigma_{2}}\,\mtx{S}^{\tau}$, exploiting that all factors are
of low rank. Then perform a low-rank update to a block-diagonal matrix whose
blocks are provided in HBS-form to construct the new HBS matrix $\mtx{T}^{\tau}$.

Accelerating the solve stage in Algorithm 2 is trivial, simply exploit that
the matrix $\mtx{S}^{\tau}$ on line (3) has low numerical rank.

\begin{remark}
Some of the structured matrix operators (e.g.~adding two HBS matrices, or the low-rank update)
can algebraically lead to a large increase in the HBS ranks. We know for
physical reasons that the output should have rank-structure very similar
to the input, however, and we combat the rank-increase by frequently
recompressing the output of each operation.
\end{remark}

\begin{remark}
\label{re:HBS}
In practice, we implement the algorithm to use dense matrix algebra at
the finer levels, where all the DtN matrices $\mtx{T}^{\tau}$ are small.
Once they get large enough that HBS algebra outperforms dense operations,
we compress the dense matrices by brute force, and rely on HBS algebra in
the remainder of the upwards sweep.
\end{remark}

\section{Numerical Examples}
\label{sec:numerics}
In this section, we illustrate the capabilities of the proposed method
for the boundary value problem

\begin{equation}
\label{eq:example}
\left\{\begin{aligned}
-\Delta u(\pxx) - c_1(\pxx)\,\partial_{1}u(\pxx)- c_2(\pxx)\,\partial_{2}u(\pxx)- c(\pxx)\,u(\pxx) =&\ 0,\qquad &\pxx \in \Omega,\\
                                     u(\pxx) =&\ f(\pxx),\qquad &\pxx \in \Gamma,
\end{aligned}\right.
\end{equation}
where $\Omega = [0,1]^2$, $\Gamma = \partial\Omega$, and $c_1(\pxx),\ c_2(\pxx),$ and $c(\pxx)$ are smooth,
cf.~(\ref{eq:defA}).
The choice of the functions $c_1(\pxx),\ c_2(\pxx),$ and $c(\pxx)$ will vary for each example.

All experiments reported in this section were executed on a machine with two quad-core Intel
Xeon E5-2643 processors and 128GB of RAM. The direct solver was implemented in Matlab, which means
that the speeds reported can very likely be improved, although the asymptotic complexity should be
unaffected.

In Section \ref{sec:known} we apply the direct solver to four problems with
known analytic solutions. This allows us to very accurately investigate the
errors incurred, but is arguably a particularly favorable environment.
Section \ref{sec:unknown} presents results from more general situations where
the exact solution is not known, and errors have to be estimated.

In all experiments, the number of Gaussian points per leaf edge $\ngauss$ is fixed at $21$,
and $21 \times 21$ tensor product grids of Chebyshev points are used in the leaf computations.
Per Remark \ref{re:HBS}, we switch from dense computations to HBS when a box has more than
$2000$ points on the boundary.

\subsection{Performance for problems with known solutions}
\label{sec:known}
To illustrate the scaling and accuracy of the discretization technique, we apply the numerical method
to four problems with known solutions.  The problems are:

\begin{description}
\item[\textit{Laplace}] Let $c_1(\pxx) = c_2(\pxx)= c(\pxx) = 0$ in (\ref{eq:example}).
\item[\textit{Helmholtz I}] Let $c_1(\pxx)= c_2(\pxx)=0$, and $c(\pxx) = \kappa^2$ where $\kappa = 80$ in (\ref{eq:example}).
This represents a vibration problem on a domain $\Omega$ of size roughly $12 \times 12$ wave-lengths. (Recall
that the wave-length is given by $\lambda = \frac{2\pi}{\kappa}$.)
\item[\textit{Helmholtz II}] Let $c_1(\pxx)= c_2(\pxx)=0$, and $c(\pxx) = \kappa^2$ where $\kappa = 640$ in (\ref{eq:example}).
This corresponds to a domain of size roughly $102 \times 102$ wave-lengths.
\item[\textit{Helmholtz III}] We again set $c_1(\pxx)= c_2(\pxx)=0$, and $c(\pxx) = \kappa^2$ in (\ref{eq:example}),
but now we let $\kappa$ grow as the number of discretization points grows to maintain a constant $12$ points per wavelength.
\end{description}

%

The boundary data in (\ref{eq:example}) is chosen to coincide with the known solutions
$u_{\rm exact}(\pxx) = \log|\hat{\pxx}-\pxx|$ for the Laplace problem and with
$u_{\rm exact}(\pxx) = \mathrm{Y}_0(\kappa|\hat{\pxx}-\pxx|)$ for the three Helmholtz
problems, where $\hat{\pxx} =(-2,0)$, and where $\mathrm{Y}_{0}$ denotes the 0'th Bessel function of the second kind.

In a first experiment, we prescribed the tolerance in the ``fast'' matrix algebra to be $\epsilon = 10^{-7}$.
Table \ref{tab:known} reports the following quantities:

\begin{tabular}{ll}
$N$              & Number of Gaussian discretization points. \\[0.5mm]
$N_{\rm tot}$    & Total number of discretization points. ($N$ plus the number of Chebyshev points)\\[0.5mm]
$T_{\rm build}$  & Time for building the solution operator.\\[0.5mm]
$T_{\rm solve} $ & Time to solve for interior nodes.\\[0.5mm]
$T_{\rm apply} $ & Time to apply the approximate Dirichlet-to-Neumann operator $\TT^1$.\\[0.5mm]
$R$              & Amount of memory required to store the solution operator.\\
$E_{\rm pot}  =$& $\max_{k\colon\pxx_{k} \in \Omega} \big\{\bigl|u_{\rm app}(\pxx_{k}) - u_{\rm exact}(\pxx_{k})\bigr|\big\},$
\end{tabular}

where $u_{\rm app}$ denotes the approximate solution constructed by the direct solver.

\begin{table}[ht]
\small
\begin{tabular}{|l|l|l|l|l|l|l|l|}
\hline
		     & $N_{\rm tot}$ &$N$      & $T_{\rm build}$ & $T_{\rm solve} $&$T_{\rm apply}$ &  $R$ & $E_{\rm pot}$\\
& && (seconds) & (seconds) & (seconds) &(MB) &\\ \hline
\multirow{3}{*}{\textit{Laplace}} &1815681 &174720   &91.68            &0.34       &0.035       &  1611.19& 2.57e-05\\
			     &7252225 &693504   &371.15		  &1.803    &0.104         &	6557.27 & 6.55e-05\\
                             & 28987905&2763264  &1661.97          &6.97  &0.168	     &	26503.29 & 2.25e-04\\
                             &115909633 &11031552 &6894.31	  &30.67   &0.367          &  106731.61 & 8.62e-04\\ \hline
\multirow{3}{*}{\textit{Helmholtz I}}&1815681  &174720   &62.07		  &0.202     &0.027        &  1611.41& 6.86e-06 \\
			     &7252225&693504   &363.19		  &1.755      &   0.084    &  6557.12& 7.47e-06\\
			      & 28987905 & 2763264 &1677.92	  &6.92      &0.186        &  26503.41 & 1.55e-05\\
                            &115909633   &11031552 &7584.65	  & 31.85    &0.435        &  106738.85& 1.45e-04\\ \hline
\multirow{3}{*}{\textit{Helmholtz II}}&1815681 &174720   &93.96		  &0.29    &0.039          &  1827.72& 5.76e-07 \\
			     &7252225 &693504   &525.92           &2.13      &0.074        &  7151.60 & 7.06e-07\\
			     & 28987905  &2763264  &2033.91	  &8.59       &0.175       &  27985.41& 4.04e-06\\ \hline
\multirow{3}{*}{\textit{Helmholtz III}}&1815681&174720   &93.68          &0.29  &    0.038        &  1839.71 & 1.29e-06\\
                            &7252225  &693504   &624.24          &1.67       &  0.086   &  7865.13 & 1.21e-06 \\
                           & 28987905    &2763264  &4174.91          &10.28   &   0.206       &  33366.45 & 1.76e-06\\ \hline
\end{tabular}
\vspace{1mm}
\caption{\label{tab:known} Timing results in seconds for the PDEs considered in Section \ref{sec:known}. For these experiments, $\epsilon = 10^{-7}$.}
\end{table}

\begin{figure}[ht]
\centering
\setlength{\unitlength}{1mm}
\begin{picture}(150,150)
\put(-10,75){\includegraphics[width=90mm]{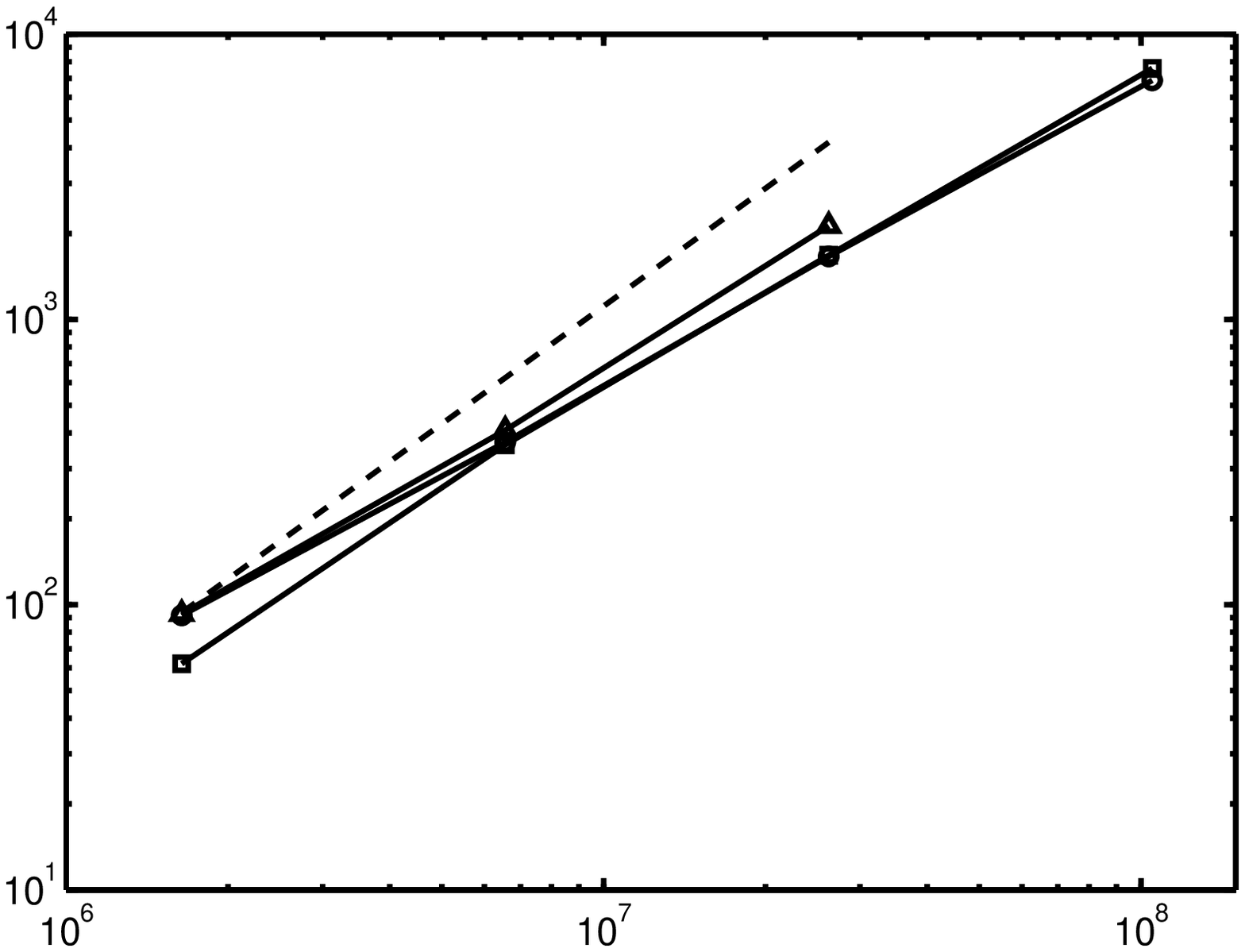}}
\put(75,75){\includegraphics[width=90mm]{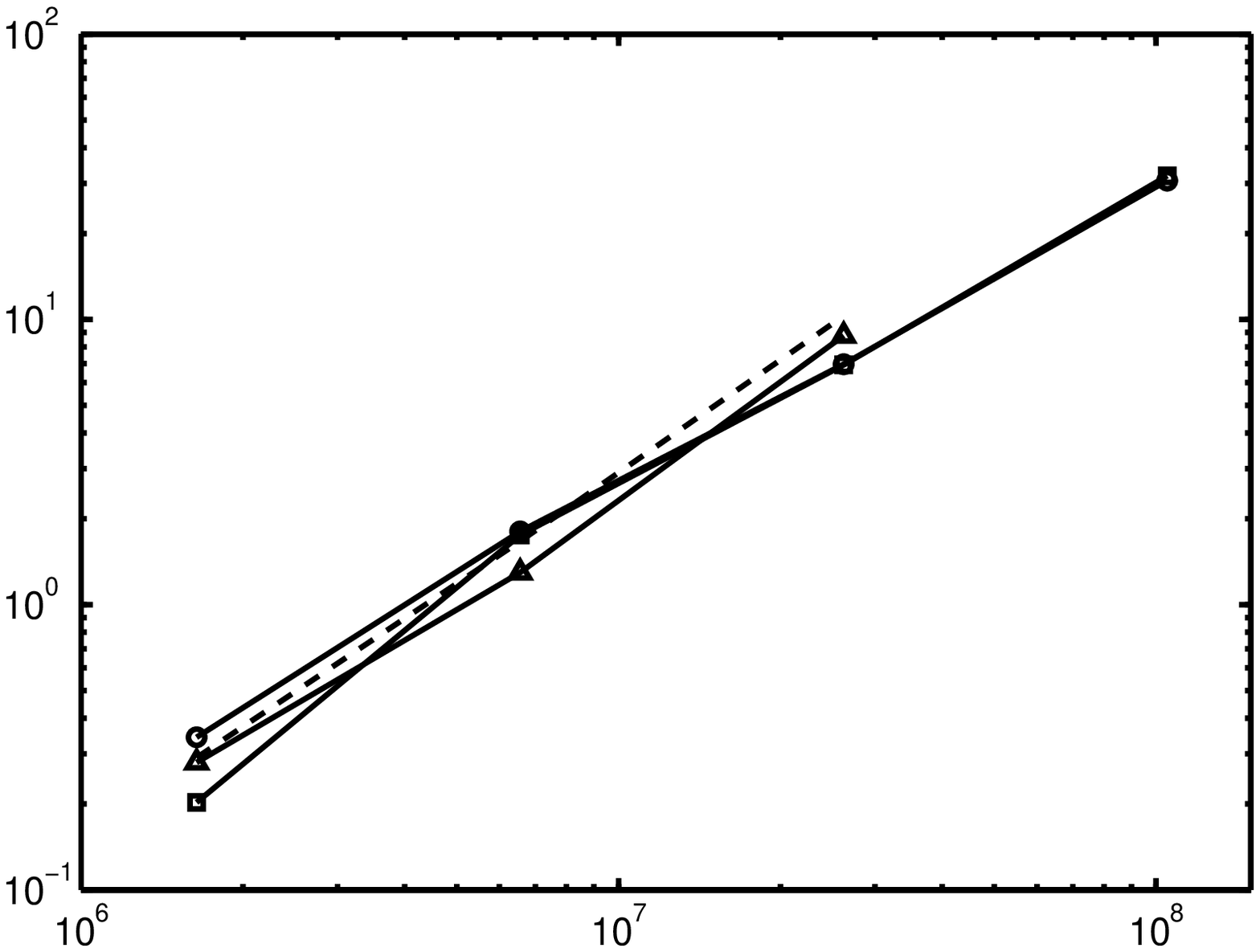}}
\put(-10,05){\includegraphics[width=90mm]{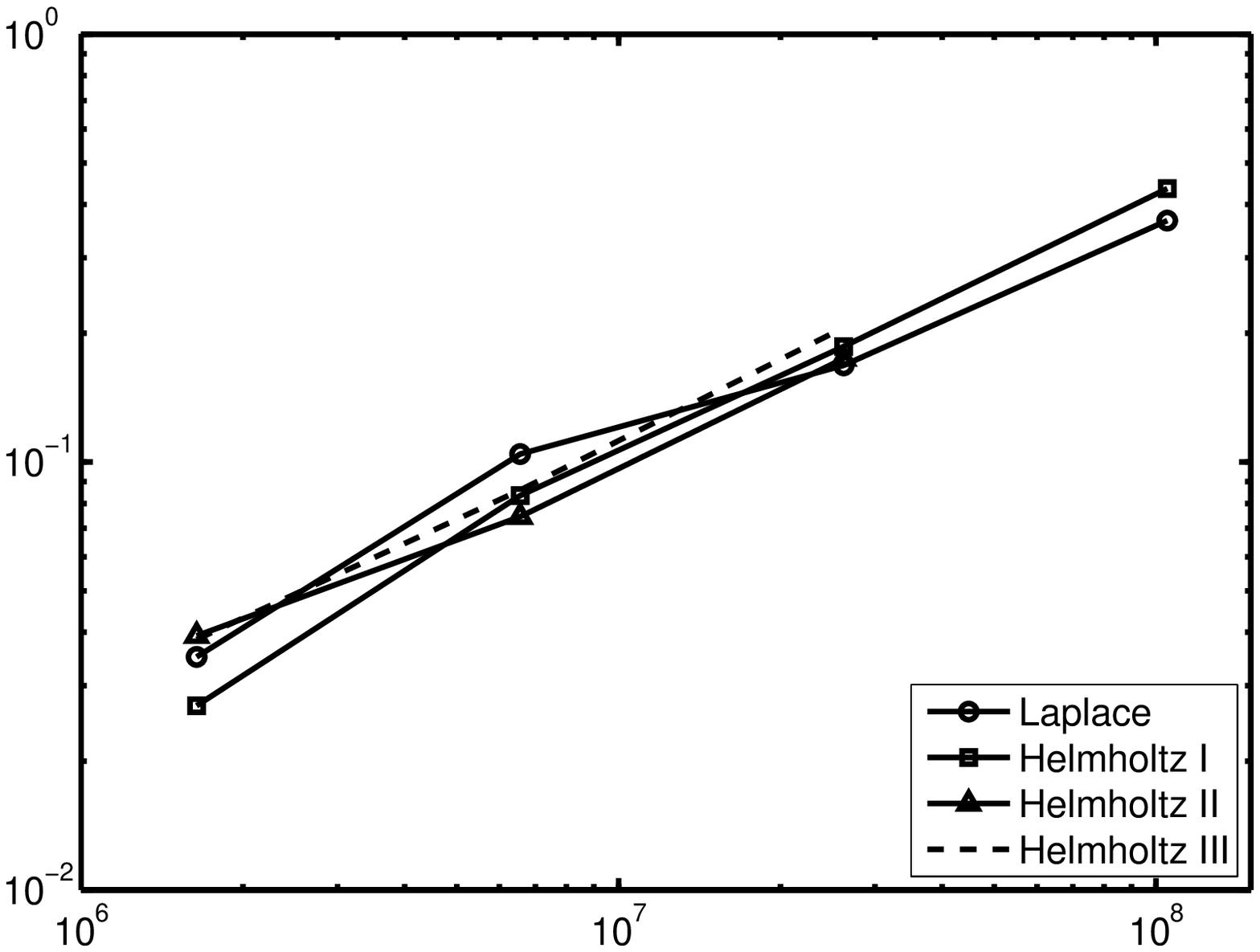}}
\put(75,05){\includegraphics[width=90mm]{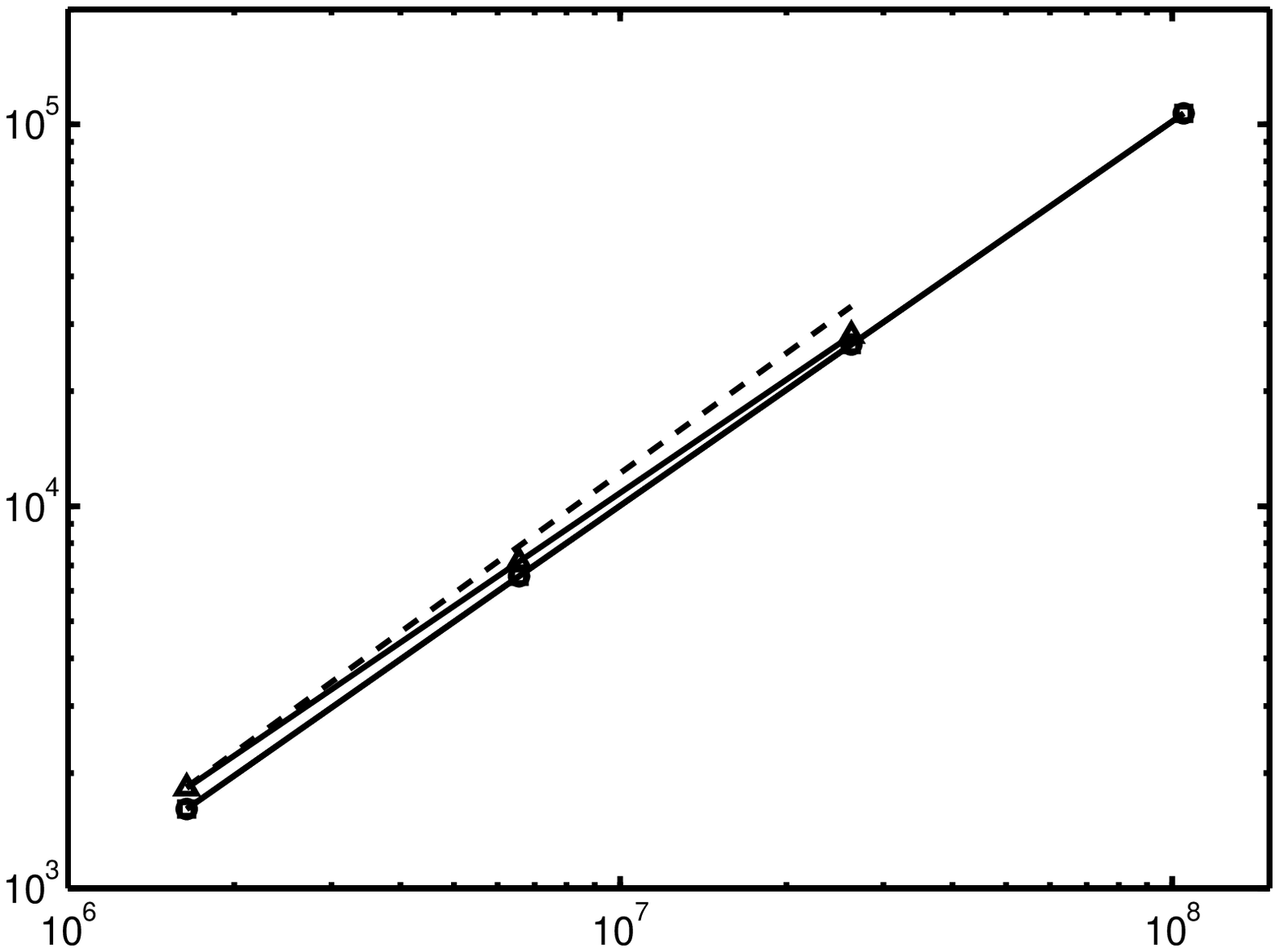}}
\put(35,135){$T_{\rm build}$}
\put(114,135){$T_{\rm solve}$}
\put(35,65){$T_{\rm apply}$}
\put(114,65){$R$ }
\put(35,05){$N$}
\put(-08,25){\rotatebox{90}{\footnotesize Time in seconds}}
\put(77,25){\rotatebox{90}{\footnotesize Memory in MB}}
\end{picture}
\caption{The first three graphs give the times required for building the direct solver ($T_{\rm build}$),
solving a problem ($T_{\rm solve}$) and applying the approximate Dirichlet-to-Neumann operator on $\partial \Omega$ ($T_{\rm apply}$).
The fourth graph gives the memory $R$ in MB required to store the solution operator.}
\label{fig:times}
\end{figure}

Our expectation is for all problems except \textit{Helmholtz III} to exhibit optimal linear
scaling for both the build and the solve stages. Additionally, we expect the cost of
applying the Dirichlet-to-Neumann operator $\mtx{T}^1$ to scale as $N^{0.5}$ for all
problems except \textit{Helmholtz III}.  The numerical results clearly bear
this out for \textit{Laplace} and \textit{Helmholtz I}. For \textit{Helmholtz II}, it
appears that linear scaling has not quite taken hold for the range of problem sizes our
hardware could manage. The \textit{Helmholtz III} problem clearly does not exhibit linear
scaling, but has not yet reached its $O(N^{1.5})$ asymptotic at the largest problem considered, which
was of size roughly $426 \times 426$ wave-lengths. We remark that the cost of the solve stage
is tiny. For example, a problem involving $11$ million unknowns (corresponding to approximately $100$
million discretization points) takes $115$ minutes for the build stage and then only $30$
seconds for each additional solve.  The cost for applying the Dirichlet-to-Neumann operator
is even less at $0.36$ seconds per boundary condition.  Figure \ref{fig:times} illustrates
the scaling via log-log plots.


In a second set of experiments, we investigated the accuracy of the computed solutions, and in
particular how the accuracy depends on the tolerance $\epsilon$ in the fast matrix algebra.
In addition to reporting $E_{\rm pot}$, Table \ref{tab:acc} reports
$$E_{\rm grad} = \max_{k\colon\pxx_{k} \in \Gamma}
\bigl\{\bigl|u_{n,\rm app}(\pxx_{k}) - u_{n,\rm exact}(\pxx_{k})\bigr|\bigr\},
$$ where $u_{\rm app}$ denotes the approximate solution constructed by the direct solver
for tolerances $\varepsilon = 10^{-7},\, 10^{-10},$ and $10^{-12}$. The number of
discretization points was fixed problem to be $N = 693504$ ($N_{\rm tot} = 7252225$).

\begin{table}[ht]
\small
\begin{tabular}{|l|l|l|l|l|l|l|}
\hline
             &  \multicolumn{2}{|c|}{$\epsilon = 10^{-7}$}& \multicolumn{2}{|c|}{$\epsilon = 10^{-10}$}& \multicolumn{2}{|c|}{$\epsilon = 10^{-12}$}\\
             &$E_{\rm pot}$ & $E_{\rm grad}$&   $E_{\rm pot}$ & $E_{\rm grad}$ &$E_{\rm pot}$ & $E_{\rm grad}$\\ \hline
\textit{Laplace}      &  6.55e-05    & 1.07e-03      & 2.91e-08        & 5.52e-07       &  1.36e-10   & 8.07e-09 \\  \hline
\textit{Helmholtz I}  &  7.47e-06    & 6.56e-04      & 5.06e-09        & 4.89e-07       & 1.38e-10    & 8.21e-09 \\ \hline
\textit{Helmholtz II} & 7.06e-07     & 3.27e-04      &  1.42e-09       & 8.01e-07       & 8.59e-11    & 4.12e-08 \\ \hline
\textit{Helmholtz III}& 1.21e-06    & 1.28e-03      &1.85e-09         & 2.69e-06       & 1.63e-09    & 2.25e-06\\ \hline
\end{tabular}

\vspace{1mm}

\caption{\label{tab:acc} Errors for solving the PDEs in Section \ref{sec:known} with different user prescribed tolerances when
the number of discretization points is fixed at $N = 693504$ ($N_{\rm tot} = 7252225$).}
\end{table}

The solution obtains (or nearly obtains) the prescribed tolerance while the normal
derivative suffers from a three digit loss in accuracy.  This loss is likely
attributable to the unboundedness of the Dirichlet-to-Neumann map.  The compressed
representation captures the high end of the spectrum to the desired accuracy while the
low end of the spectrum is captured to three digits less than the desired accuracy.

\subsection{Convergence for unknown solutions}
\label{sec:unknown}
In this section, we apply the direct solver to three problems for which we
do not know an exact solution:

\begin{description}
\item[\textit{Constant convection}]
Let the convection in the $x_2$ direction be constant by setting $c_1(\pxx) = 0$,
$c_2(\pxx) = -100,000$ and set $c(\pxx) = 0$.
\item[\textit{Diffusion-Convection}]
Introduce a divergence free convection by setting
$c_1(\pxx) = -10,000\cos(4 \pi x_2)$, $c_2(\pxx) = -10,000\cos(4 \pi x_1)$, and $c(\pxx) = 0$.
\item[\textit{Variable Helmholtz}] Consider the variable coefficient Helmholtz problem where
$c_1(\pxx) = 0$, $c_2(\pxx) = 0$, $c(\pxx) = \kappa^2(1-(\sin(4\pi x_1)\sin(4 \pi x_2))^2)$
and $\kappa = 640$.
\end{description}

\noindent
For the three experiments, the boundary data is given by $f(\pxx) = \cos(2 x_1)(1-2x_2)$.

To check for convergence, we post-process the solution as described in Section \ref{sec:postproc}
to get the solution on the Chebyshev grid.  Let $u^N$ denote the solution on the
Chebyshev grid.  Likewise,let $u_n^N$ denote the normal derivative on the boundary
at the Chebyshev boundary points.  We compare the solution and the normal
derivative on the boundary pointwise at the locations
$$ \hat{\pxx} = (0.75,0.25)\quad {\rm and}\quad \hat{\pyy} = (0.75,0) $$
respectively, via
$$ E^N_{\rm int} = |u^N(\hat{\pxx})-u^{4N}(\hat{\pxx})| \quad {\rm and} \quad E^N_{\rm bnd} = |u_n^N(\hat{\pyy})-u_n^{4N}(\hat{\pyy})|.$$

The tolerance for the compressed representations is set to $\epsilon = 10^{-12}$.
Table \ref{tab:conv} reports the pointwise errors. We see that high accuracy is
obtained in all cases, with solutions that have ten correct digits for the potential
and about seven correct digits for the boundary flux.

The computational costs of the computations are reported in Table \ref{tab:times}.
The memory $R$ reported now includes the memory required to store all local
solution operators described in Section \ref{sec:postproc}.

\begin{table}
 \small
\begin{tabular}{|l|l|l|l|l|l|l|}
\hline
                              & $N_{\rm tot}$   & $N$   & $u^{N}(\hat{\pxx})$ & $E^{N}_{\rm int}$& $u_n^{N}(\hat{\pxx})$ & $E^{N}_{\rm bnd}$\\ \hline
\multirow{4}{*}{\textit{Constant Convection}} & 455233& 44352     & -0.474057246846306     & 0.477            &-192794.835134257        &824.14       \\
                                     &1815681 &174720 & -0.951426960146812      &  8.28e-03          &-191970.688228300        &1.47 \\
                                     &7252225  &693504  & -0.959709514830931     &  6.75e-10          & -191972.166912008        &0.365 \\
                                     & 28987905 &2763264& -0.959709515505929     &                    & -191972.532606428       &   \\ \hline
\multirow{4}{*}{\textit{Variable Helmholtz}}&114465 &11424    & 2.50679456864385      &  6.10e-02          &-2779.09822864819        & 3188 \\
                                    & 455233& 44352   &2.56780367343056       &  4.63e-07          & 409.387483435691        & 2.59e-02  \\
                                     &1815681 &174720  &2.56734097240752       &  1.77e-09          &409.413356177218         & 3.30e-07 \\
                                   &7252225  &693504  & 2.56734097418159      &                    &409.413355846946         & \\ \hline
\multirow{4}{*}{\textit{Diffusion-Convection}}& 455233& 44352   & 0.0822281612709325    &5.04e-5             &-35.1309711271060        & 2.23e-3          \\
				     &1815681 &174720  & 0.0822785917678385    &2.67e-8             &-35.1332056731696        &7.57e-6  \\
				     &7252225  &693504  & 0.0822786184699243    &5.41e-12            &-35.1332132455725        & 2.11e-09 \\
				      & 28987905 &2763264& 0.0822786184753420    &                    &-35.1332132476795	      & \\ \hline
\end{tabular}
\vspace{1mm}
\caption{\label{tab:conv}Convergence results for solving the PDE's in Section \ref{sec:unknown} with a user subscribed tolerance of $\epsilon = 10^{-12}$.}
\end{table}

\begin{table}
\small
\begin{tabular}{|l|l|l|l|l|l|}
\hline
                             & $N_{\rm tot}$    & $N$    &$T_{\rm build}$ & $T_{\rm solve}$ & $R$\\
&&& (seconds) & (seconds) & (MB) \\ \hline
\multirow{4}{*}{\textit{Constant Convection}}& 455233& 44352     &21.04          & 0.85           & 683.25\\
				    &1815681 &174720     &176.09          &3.47             &2997.80\\
				    &7252225  &693504     &980.93         &13.76            &8460.94\\
				    & 28987905 &2763264  &5227.52         & 77.03           & 48576.75\\ \hline
\multirow{4}{*}{\textit{Variable Helmholtz}} &114465 &11424      &4.61            &0.19             & 167.68 \\
                                    & 455233& 44352     &42.72           & 1.110           & 774.34\\
                                     &1815681 &174720  &450.68          &4.54             & 3678.31\\
                                     &7252225  &693504     &3116.57         &17.64            &15658.07\\ \hline
\multirow{4}{*}{\textit{Diffusion-Convection}}& 455233& 44352   &28.31           & 0.795           & 446.21\\
				     &1815681 &174720   &131.23          &3.20             & 2050.20\\
				     &7252225  &693504    &906.11          &17.12            &8460.94\\
				     & 28987905 &2763264 &4524.99         & 66.99           & 47711.17 \\ \hline
\end{tabular}

\vspace{1mm}
\caption{\label{tab:times}Times in seconds for solving the PDE's in Section \ref{sec:unknown} with a user subscribed tolerance of $\epsilon = 10^{-12}$.}
\end{table}


\section{Conclusions}
\label{sec:conc}
We have described a direct solver for variable coefficient elliptic PDEs in the plane,
under the assumption that the solution and all coefficient functions are smooth.
For problems with non-oscillatory solutions such as the Laplace and Stokes equations,
the asymptotic complexity of the solver is $O(N)$, with a small constant of proportionality.
For problems with oscillatory solutions, high practical efficiency is retained for
problems of size up to several hundred wave-lengths.

Our method is based on a composite spectral discretization. We use high order local
meshes (typically of size $21\times 21$) capable of solving even very large scale
problems to ten correct digits or more. The direct solver is conceptually similar
to the classical nested dissection method \cite{george_1973}. To improve the asymptotic
complexity from the classical $O(N^{1.5})$ to $O(N)$, we exploit internal structure in
the dense ``frontal matrices'' at the top levels in a manner similar to recent work such
as, e.g., \cite{2010_jianlin_fast_hss,Adiss,2007_leborne_HLU,2009_martinsson_FEM,2011_ying_nested_dissection_2D}.
Compared to these techniques, our method has an advantage in that high order discretizations
can be incorporated without compromising the speed of the linear algebra. The reason is
that we use a formulation based on Dirichlet-to-Neumann operators. As a result, we need
high order convergence \textit{only in the tangential direction} on patch boundaries.

The direct solver presented requires more storage than classical iterative methods,
but this is partially off-set by the use of high-order discretizations. More
importantly, the solver is characterized by very low data-movement. This makes the
method particularly well suited for implementation on parallel machines with distributed
memory.

\vspace{5mm}

\noindent\textbf{Acknowledgments:}
The work was supported by NSF grants DMS-0748488 and CDI-0941476.


\bibliography{main_bib}

\providecommand{\bysame}{\leavevmode\hbox to3em{\hrulefill}\thinspace}
\providecommand{\MR}{\relax\ifhmode\unskip\space\fi MR }
\providecommand{\MRhref}[2]{%
  \href{http://www.ams.org/mathscinet-getitem?mr=#1}{#2}
}
\providecommand{\href}[2]{#2}
\begin{thebibliography}{10}

\bibitem{2004_gu_divide}
S.~Chandrasekaran and M.~Gu, \emph{A divide-and-conquer algorithm for the
  eigendecomposition of symmetric block-diagonal plus semiseparable matrices},
  Numer. Math. \textbf{96} (2004), no.~4, 723--731.

\bibitem{2010_jianlin_fast_hss}
S.~Chandrasekaran, M.~Gu, X.S. Li, and J~Xia, \emph{Fast algorithms for
  hierarchically semiseparable matrices}, Numer. Linear Algebra Appl.
  \textbf{17} (2010), 953--976.

\bibitem{2013_yu_chen_totalwave}
Y.~Chen, \emph{Total wave based fast direct solver for volume scattering
  problems}, ar{X}iv.org \textbf{1302.2101} (2013).

\bibitem{2002_chen_direct_lippman_schwinger}
Yu~Chen, \emph{A fast, direct algorithm for the {L}ippmann-{S}chwinger integral
  equation in two dimensions}, Adv. Comput. Math. \textbf{16} (2002), no.~2-3,
  175--190, Modeling and computation in optics and electromagnetics.

\bibitem{1989_directbook_duff}
I.S. Duff, A.M. Erisman, and J.K. Reid, \emph{Direct methods for sparse
  matrices}, Oxford, 1989.

\bibitem{george_1973}
A.~George, \emph{Nested dissection of a regular finite element mesh}, SIAM J.
  on Numerical Analysis \textbf{10} (1973), 345--363.

\bibitem{Adiss}
A.~Gillman, \emph{Fast direct solvers for elliptic partial differential
  equations}, Ph.D. thesis, University of Colorado at Boulder, Applied
  Mathematics, 2011.

\bibitem{m2011_1D_survey}
A.~Gillman, P.~Young, and P.G. Martinsson, \emph{A direct solver with $o(n)$
  complexity for integral equations on one-dimensional domains}, Frontiers of
  Mathematics in China \textbf{7} (2012), no.~2, 217--247.

\bibitem{2009_martinsson_ACTA}
Leslie Greengard, Denis Gueyffier, Per-Gunnar Martinsson, and Vladimir Rokhlin,
  \emph{Fast direct solvers for integral equations in complex three-dimensional
  domains}, Acta Numer. \textbf{18} (2009), 243--275.

\bibitem{2012_ho_greengard_fastdirect}
K.L. Ho and L.~Greengard, \emph{A fast direct solver for structured linear
  systems by recursive skeletonization}, SIAM J Sci Comput \textbf{to appear}
  (2012).

\bibitem{hoffman_1973}
A.~J. Hoffman, M.~S. Martin, and D.~J. Rose, \emph{Complexity bounds for
  regular finite difference and finite element grids}, SIAM J. Numer. Anal.
  \textbf{10} (1973), 364--369.

\bibitem{1998_kopriva_multidomain}
D.A. Kopriva, \emph{A staggered-grid multidomain spectral method for the
  compressible navier–stokes equations}, Journal of Computational Physics
  \textbf{143} (1998), no.~1, 125 -- 158.

\bibitem{2007_leborne_HLU}
Sabine Le~Borne, Lars Grasedyck, and Ronald Kriemann,
  \emph{Domain-decomposition based {$\mathcal{H}$}-{LU} preconditioners},
  Domain decomposition methods in science and engineering {XVI}, Lect. Notes
  Comput. Sci. Eng., vol.~55, Springer, Berlin, 2007, pp.~667--674.
  \MR{2334161}

\bibitem{2009_martinsson_FEM}
P.G. Martinsson, \emph{A fast direct solver for a class of elliptic partial
  differential equations}, J. Sci. Comput. \textbf{38} (2009), no.~3, 316--330.

\bibitem{2012_martinsson_spectralcomposite}
\bysame, \emph{A composite spectral scheme for variable coefficient helmholtz
  problems}, ar{X}iv.org \textbf{1206.4136} (2012).

\bibitem{2012_martinsson_spectralcomposite_jcp}
\bysame, \emph{A direct solver for variable coefficient elliptic pdes
  discretized via a composite spectral collocation method}, Journal of
  Computational Physics \textbf{242} (2013), no.~1, 460 -- 479.

\bibitem{2005_martinsson_fastdirect}
P.G. Martinsson and V.~Rokhlin, \emph{A fast direct solver for boundary
  integral equations in two dimensions}, J. Comp. Phys. \textbf{205} (2005),
  no.~1, 1--23.

\bibitem{2003_pfeiffer_spectralmultidomain}
H.P. Pfeiffer, L.E. Kidder, M.A. Scheel, and S.A. Teukolsky, \emph{A
  multidomain spectral method for solving elliptic equations}, Computer physics
  communications \textbf{152} (2003), no.~3, 253--273.

\bibitem{2011_ying_nested_dissection_2D}
P.G. Schmitz and L.~Ying, \emph{A fast direct solver for elliptic problems on
  general meshes in 2d}, Journal of Computational Physics \textbf{231} (2012),
  no.~4, 1314 -- 1338.

\bibitem{2007_shiv_sheng}
Zhifeng Sheng, Patrick Dewilde, and Shivkumar Chandrasekaran, \emph{Algorithms
  to solve hierarchically semi-separable systems}, System theory, the Schur
  algorithm and multidimensional analysis, Oper. Theory Adv. Appl., vol. 176,
  Birkh\"auser, Basel, 2007, pp.~255--294. \MR{MR2342902}

\bibitem{2000_trefethen_spectral_matlab}
L.N. Trefethen, \emph{Spectral methods in matlab}, SIAM, Philadelphia, 2000.

\bibitem{2009_xia_superfast}
Jianlin Xia, Shivkumar Chandrasekaran, Ming Gu, and Xiaoye~S. Li,
  \emph{Superfast multifrontal method for large structured linear systems of
  equations}, SIAM J. Matrix Anal. Appl. \textbf{31} (2009), no.~3, 1382--1411.
  \MR{2587783 (2011c:65072)}

\bibitem{2000_hesthaven_pseudospectral}
B.~Yang and J.S. Hesthaven, \emph{Multidomain pseudospectral computation of
  maxwell's equations in 3-d general curvilinear coordinates}, Applied
  Numerical Mathematics \textbf{33} (2000), no.~1 -- 4, 281 -- 289.

\end{thebibliography}
\bibliographystyle{amsplain}

\end{document}